\documentclass[12pt]{article}

\usepackage{amssymb,amsmath} 

\numberwithin{equation}{section}

\usepackage{amsfonts}

\usepackage{enumerate} 

\usepackage{color}
\definecolor{lightblue}{rgb}{0,0.2,0.5}
\definecolor{newcolor}{rgb}{0,0,0}

\usepackage[colorlinks=true, urlcolor=blue,linkcolor=blue, citecolor=lightblue]{hyperref}

\usepackage{graphics,graphicx,amsmath}

\usepackage{empheq} 

\def\cF{\mathcal F}
\def\cE{\mathcal E}

\newcommand{\R}{\mathbb{R}}

\newcommand{\E}{\mathbb{E}}

\newtheorem{prop}{Proposition}[section]
\newtheorem{lemma}[prop]{Lemma}
\newtheorem{definition}[prop]{Definition}
\newtheorem{corollary}[prop]{Corollary}
\newtheorem{theorem}[prop]{Theorem}
\newtheorem{remark}[prop]{Remark}
\newtheorem{example}[prop]{Example}

\def\({\left(}
\def\){\right)}

\def\[{\left[}
\def\]{\right]}
\def\real{{\mathord{\mathbb R}}}

\makeatletter
\newcommand*\rel@kern[1]{\kern#1\dimexpr\macc@kerna}
\newcommand*\widebar[1]{
  \begingroup
  \def\mathaccent##1##2{
    \rel@kern{0.8}
    \overline{\rel@kern{-0.8}\macc@nucleus\rel@kern{0.2}}
    \rel@kern{-0.2}
  }
  \macc@depth\@ne
  \let\math@bgroup\@empty \let\math@egroup\macc@set@skewchar
  \mathsurround\z@ \frozen@everymath{\mathgroup\macc@group\relax}
  \macc@set@skewchar\relax
  \let\mathaccentV\macc@nested@a
  \macc@nested@a\relax111{#1}
  \endgroup
}
\makeatother

\usepackage{dsfont}

\usepackage[round]{natbib} 

\bibpunct{\textcolor{lightblue}{(}}{\textcolor{lightblue}{)}}{,}{a}{}{;}

\let\oldcitet=\citet
\let\oldcitep=\citep 
\renewcommand{\cite}[1]{\textcolor[rgb]{0,0,1}{\oldcitet{#1}}}
\renewcommand{\citet}[1]{\textcolor[rgb]{0,0,1}{\oldcitet{#1}}}
\renewcommand{\citep}[1]{\textcolor[rgb]{0,0,1}{\oldcitep{#1}}}

\newenvironment{Proof}{\removelastskip\par\medskip
\noindent{\em Proof.} \rm}{\penalty-20\null\hfill$\square$\par\medbreak}

\allowdisplaybreaks

\begin{document}
\title{
\Huge
Stochastic ordering by $g$-expectations 
} 

\author{Sel Ly \qquad Nicolas Privault
  \\
  \vspace{-0.1cm}
\small
School of Physical and Mathematical Sciences
\\
\vspace{-0.1cm}
\small
Division of Mathematical Sciences
\\
\vspace{-0.1cm}
\small
 Nanyang Technological University \\
\vspace{-0.1cm}
\small
Singapore 637371
}

\maketitle

\vspace{-0.7cm} 

\begin{abstract} 
 We derive sufficient conditions for the
 convex and monotonic $g$-stochastic ordering of diffusion
 processes under nonlinear $g$-expectations and $g$-evaluations.
 Our approach relies on comparison results for forward-backward 
 stochastic differential equations  
 and on several extensions of convexity, monotonicity and continuous dependence 
 properties for the solutions of associated 
 semilinear parabolic partial differential equations.
 Applications to contingent claim price comparison under
 different hedging portfolio constraints are provided. 
\end{abstract}

\baselineskip0.7cm

\section{Introduction}
In comparison with standard mean-variance analysis, 
partial orderings of probability distributions
provide additional information which can be  
used in applications to risk management. 
In this framework, a random variable $X^{(1)}$ is said 
to be dominated by another random variable $X^{(2)}$ if 
\begin{equation} 
\label{111} 
 \E \big[\phi\big(X^{(1)}\big)\big] \leq  \E \big[\phi\big(X^{(2)}\big)\big] , 
\end{equation} 
 for all $\phi : \real \to \real$ in a certain class
 of functions,
 where $\E[X]$ denotes the usual expectation of the
 random variable $X$.
For example, if $X^{(1)}$ and $X^{(2)}$ represent the lifetimes
of two devices $A$ and $B$ then the stochastic ordering
\eqref{111} for all non-decreasing and bounded functions $\phi$,
tells that the device $B$ will likely survive longer than the device $A$.
Stochastic ordering
  has found a wide range of applications in various fields such as reliability theory, economics, actuarial sciences,  operation research, risk management, biology, option evaluation, etc., see e.g. \cite{mullerbk,denuit, shaked2007stochastic,sriboonchita2009stochastic,levy2015stochastic,belzunce2015introduction,perrakis2019stochastic}.
\textcolor{newcolor}{In the von Neumann-Morgenstern
  expected utility theory, a portfolio with
  return $X^{(2)}$ dominates a portfolio with return
  $X^{(1)}$ in the increasing concave order
  if \eqref{111} holds
  for all non-decreasing concave utility functions $\phi$,
  in which case, the second portfolio would be preferred
  over the first portfolio by risk-averse investors,
  see Theorem~1.35 in \cite{sriboonchita2009stochastic}.
 Similarly, if \eqref{111} holds
  for all non-decreasing convex utility functions $\phi$,
  the second portfolio would be preferred
  over the first portfolio by risk-seeking investors,
  see Theorem~1.37
  and the notion of stochastic dominance in Theorem~2.4
  in \cite{sriboonchita2009stochastic}.
}

\medskip 

Comparison bounds in convex ordering have been established in \cite{elkjs} 
for option prices with convex payoff functions
in the continuous diffusion case, via a martingale approach
based on the classical Kolmogorov equation 
and the propagation of convexity property for Markov semigroups.
This approach has been generalized to semimartingales
in \cite{mordecki},
\cite{bergenthum},
\cite{bergenthum2}, see also
\cite{kmp},
\cite{abp},
\cite{ma-privault}. 

\medskip 

\textcolor{newcolor}{
  On the other hand, empirical experiments have shown
  that many uncertain phenomena
  cannot be fully modeled using
  the linear expectation operator $\E[ \ \! \cdot \ \! ]$,
  as in e.g. the Allais and Ellsberg paradoxes.
  Choquet's expectation has been proposed
  as an nonlinear alternative
  that relies on capacities instead of probability measures,
  see \cite{grigorova2014stochastic,grigorova2014dominance} 
  for the construction of monotonic and increasing convex stochastic orders
  and application to financial optimization.
}

\medskip 

\textcolor{newcolor}{ The nonlinear $g$-expectation
and $g$-evaluation 
$\cE_g[\xi]$ of a random variable $\xi$
have been introduced} 
by \cite{peng1997backward,peng2004nonlinear} 
	as the initial value $ Y_0$ for a pair $(Y_t,Z_t)_{t\in [0,T]}$
	of adapted processes
 solution of a Backward Stochastic Differential Equation (BSDE) of the form 
	\begin{equation}
\nonumber 
		-dY_t = g(t, X_t, Y_t, Z_t)dt - Z_tdB_t, \qquad 0 \leq t \leq T, 
	\end{equation}
	with terminal condition $Y_T = \xi$, where the function $g(t,x,y,z) $ is called the BSDE
	generator, $(B_t)_{t\in \real_+}$ is a standard Brownian motion defined on a probability measure space $\big(\Omega, \cF,\mathbb{P}\big)$, 
	and $(X_t)_{t\in \real_+}$ is a diffusion process driven by $(B_t)_{t\in \real_+}$. The $g$-expectation $\cE_g$  preserves all properties of the classical expectation $\E$, except for linearity, and it generalizes the classical notion of
	expectation which corresponds to the choice $ g(t,x,y,z) :=0$,
	$t\in [0,T]$, $x,y,z\in \real$.

\medskip 

BSDEs were first introduced by \cite{bismut1973conjugate} in
the linear case, and then extended by \cite{pardoux1990adapted}
to the nonlinear case. \textcolor{newcolor}{BSDEs and the corresponding $g$-expectations $\cE_g$ have been applied to contingent claim pricing,
stochastic control theory,
utility maximization and dynamic risk measures, see e.g. 
\cite{pardoux1990adapted,peng1997backward,el1997backward,ma1999forward, 
  peng2004nonlinear,pengicm, gianin2006risk, epstein2013ambiguous, epstein2014ambiguous, jiang2016utility}. 
}

\medskip 

\textcolor{newcolor}{
In this paper, we study stochastic orderings from the point of view of
nonlinear $g$-expectations.
Consider two risky assets with positive prices
$\big(X_t^{(i)}\big)_{t\in [0,T]}$, $i=1,2$, given by
\begin{equation}
	\nonumber
	dX_t^{(i)}=X_t^{(i)} a_i\big(t,X_t^{(i)}\big)dt + X_t^{(i)} b_i\big( t, X_t^{(i)} \big) dB_t,
	\qquad
	i=1,2, 
\end{equation}
a risk-free asset $E_t := E_0e^{rt}$, where $r$ is an interest rate,
and two portfolios with prices
$Y^{(i)}_t = p_t^{(i)} E_t + q_t^{(i)} X_t^{(i)}$,
 under the self-financing conditions 
\begin{align}
	dY_t^{(i)} &  =q^{(i)}_tdE_t + p^{(i)}_tdX_t^{(i)} \nonumber
	\\
	& =\big(
	rY_t^{(i)} +
	\big(a_i\big( t,X_t^{(i)}\big)-r\big)
	p^{(i)}_t X_t^{(i)} 
	\big) dt +
	p^{(i)}_tX_t^{(i)} b_i \big( t,X_t^{(i)}\big)dB_t, \quad i=1,2, 
	 \end{align} 
 which lead to the BSDEs 
\begin{equation}
	\nonumber
	Y_t^{(i)}=Y_T^{(i)}
	+ \int_t^T g_i\big(s,X_s^{(i)}, Y_s^{(i)}, Z_s^{(i)}\big) ds
	- \int_t^T Z_s^{(i)}dB_s,
\end{equation}
where $Z_t^{(i)}:=p^{(i)}_tX_t^{(i)} b_i \big( t,X_t^{(i)}\big)$ and
$$
g_i(t,x,y,z) := -ry- z \theta_i( t,x ),
\quad
\mbox{with}
\quad
\theta_i (t,x) = \frac{a_i( t,x )-r}{b_i( t,x)},
\quad i=1,2.
$$ 
 We say that $X_T^{(2)}$ dominates $X_T^{(1)}$ in the convex $g_1,g_2$-stochastic ordering,
i.e. $X_T^{(1)} \le_{g_1,g_2}^{conv} X_T^{(2)}$ if
$$\cE_{g_1}\big[\phi\big(X_T^{(1)}\big)\big] \leq \cE_{g_2}\big[\phi\big(X_T^{(2)}\big)\big]
$$
for all convex functions $\phi$,
where  $Y_0^{(i)}:=\cE_{g_i}[\phi\big(X_T^{(i)}\big)]$
represent the fair prices at time $t=0$ of the options with convex payoffs 
$Y_T^{(i)}:=\phi\big(X_T^{(i)}\big)$, $i=1,2$. }

\medskip

\textcolor{newcolor}{
 Here, the use of distinct generators
 $g_i(t,x,y,z)$ is motivated by the comparison
  of different contingent claims
    under different hedging strategies,
    for example in the case of misspecified volatility coefficients
    or for hedging under constraints,
    see the examples presented in Section~\ref{s4}.
  In case $g(t,x,y,z)=\beta(t)z$,
  the increasing convex $g$-stochastic ordering $\le^{iconv}_g$ is
  equivalent to the (classical) increasing convex ordering with respect
  to the capacities $\mu_g$,
  see \cite{grigorova2014stochastic}, where $\mu_g[A]:=\cE_g[\mathbf{1}_A]$,
  $A \in \cF$, see \cite{chen-davison}.}

\subsubsection*{Main results} 

  In Theorem~\ref{pp3}, 
        we derive sufficient conditions
        on two BSDE generators
        $g_1(t,x,y,z )$, $g_2(t,x,y,z )$ 
        for the convex ordering 
        \begin{equation}
\nonumber 
  	\cE_{g_1}\big[\phi\big(X_T^{(1)}\big)\big] \leq \cE_{g_2}\big[\phi\big(X_T^{(2)}\big)\big],
  	\end{equation} 
  	in nonlinear expectations $\cE_{g_1}$, $\cE_{g_2}$,  
        for all convex functions $\phi(x)$
        with polynomial growth,
        where $X_T^{(1)}$ and $X_T^{(2)}$
        are the terminal values
        of the solutions of two 
forward Stochastic Differential Equations (SDEs)
        \begin{subequations}
        	\begin{empheq}[left=\empheqlbrace]{align}
        	\nonumber 
        		dX_t^{(1)}=\mu_1\big(t,X_t^{(1)}\big)dt + \sigma_1\big(t,X_t^{(1)}\big)dB_t, 
        		\\
        		\nonumber 
        		\\ 
\nonumber 
        		dX_t^{(2)}=\mu_2\big(t,X_t^{(2)}\big)dt + \sigma_2\big(t,X_t^{(2)}\big)dB_t, 
        	\end{empheq}
        \end{subequations}
 with $X_0^{(1)}= X_0^{(2)}$, under the bound 
  	$$
  	0 < \sigma_1(t,x) \leq \sigma_2(t,x), \qquad t\in [0,T], 
\quad x \in \real. 
  	$$
   
   \textcolor{newcolor}{
     \noindent
     The proof of Theorem~\ref{pp3} 
 is based on}
   the comparison Theorem~2.4 in Appendix~C of \cite{peng},
        provided that 
        $$
        z\mu_1(t,x)
  	+
  	g_1(t,x,y,z \sigma_1(t,x))
  	\leq
  	z\mu_2(t,x)
  	+
  	g_2(t,x,y,z\sigma_2(t,x)), 
  	\quad x,y,z \in \real, \ t\in [0,T],
  	$$
        and both functions
        $(x,y,z) \mapsto f_i(t,x,y,z):=z \mu_i(t,x)+ g_i(t,x,y,z\sigma_i(t,x))$
        are convex in $(x,y)$ and in $(y,z)$
        on $\real^2$ for $i=1,2$ and $t \in [0,T]$.       
        \textcolor{newcolor}{
          Several extensions are considered
          on increasing convex and monotonic orderings
          in Theorem~\ref{pp4} and Corollaries~\ref{pp4.1}-\ref{pp5},
            with the particular cases of equal drifts and equal
volatilities treated in Corollaries~\ref{pp6} and \ref{pp7}.}

\medskip

\textcolor{newcolor}{This approach requires convexity of the function $(x,y,z) \mapsto f_i(t,x,y,z)$ for both $i=1$ and $ i=2 $.
  In Section~\ref{s3.1} we relax those conditions
  using a stochastic calculus approach, 
  by only requiring the convexity of the function $ f_i(t,x,y,z) $ for $ i=1 $
  or $ i=2 $ in Theorems~\ref{pp1} and \ref{pp2}, 
  which respectively deal with the convex and increasing convex orders. } 
 
\medskip 

Related comparison results for $g$-risk measures are 
\textcolor{newcolor}{presented in 
Corollaries~\ref{pp9}-\ref{pp10}, 
 using the quantity
 $ \cE_{g_i}\big[-\phi\big(X_T^{(i)}\big)~\big| \cF_t\big]=-\cE_{g_i^{(-1)}}\big[\phi\big(X_T^{(i)}\big)~\big| \cF_t\big] $
 which makes sense as a dynamic $g$-risk measure, 
 where $ g_i^{(-1)}(t,x,y,z):=-g_i(t,x,-y,-z) $. 
 Here, the choice of generator function $ g_i $
 determines the investor's portfolio
 strategy and the corresponding risk measures,
 see Section~\ref{s4} for examples. 
}
 
\medskip
    
  The proofs of
    Theorems~\ref{pp3}-\ref{pp10}
  rely on an extension of convexity properties
  of the solutions of nonlinear parabolic Partial Differential Equations (PDEs)
  which is proved in Theorem~\ref{convex}.
  The convexity properties of solutions of nonlinear PDEs 
  have been studied by several authors,
  see e.g. Theorem~3.1 in \cite{lionsmusiela},
  Theorem~2.1 in \cite{giga}, 
  and Theorem~1.1 in \cite{bian2008convexity}, 
  see also Theorem~1 in \cite{alvarezlasry} in the
  elliptic case. 
  Those works typically require global convexity of the nonlinear
  drifts $f(t,x,y,z)$ in all state variables $(x,y,z)$,
  a condition which is too strong
  for our applications to finance in 
  Examples~\ref{ex6.1}-\ref{ex6.4} below. 
  For this reason, in Theorem~\ref{convex} 
  we extend Theorem~1.1 of \cite{bian2008convexity}
  in dimension one, by replacing
  the global convexity of the nonlinear drift $f_i(t,x,y,z)$ in $(x,y,z)$  
  with its convexity in $(x,y)$ and $(y,z)$, $i=1,2$.

\medskip

\textcolor{newcolor}{Finally,} Section~\ref{s6} deals with monotonicity properties
          and continuous dependence results
for the solutions
 Forward-Backward Stochastic
 Differential Equations (FBSDEs)
 and PDEs,
 which are used in the proofs of Theorems~\ref{pp3}-\ref{pp2} and
 Corollaries~\ref{pp4.1}-\ref{pp7}. 

 \section{Preliminaries}
 In this section, we recall some notation and background
          on FBSDEs,
          $g$-expectations,
          $g$-evaluations
          and $g$-stochastic orderings.           
          Given $T>0$,
          let $(B_t)_{t\in [0,T]}$ be a standard Brownian motion on
a probability space $\left(\Omega,{\cal F}, \mathbb{P} \right)$.
Denote by $\left(\cF_t\right)_{t\in [0,T]}$ the augmented filtration such that
$\cF_t = \sigma ( B_s, 0\leq s \leq t ) \vee {\cal N}$,
$t \in [0,T]$, where ${\cal N} $ is the collection of all $\mathbb{P}$-null sets. We also let $L^2(\Omega, \cF_t) := L^2(\Omega,\cF_t, \mathbb{P})$, 
$t \in [0, T]$. 
	\subsubsection*{Forward-Backward SDEs}
Consider a forward SDE of the form  
\begin{equation}
\nonumber 
	dX_s^{t,x}=\mu\big(s,X_s^{t,x}\big)dt+\sigma\big(s, X_s^{t,x}\big)dB_s, \qquad 0\leq t \leq s \leq T. 
\end{equation}
with initial condition $X_t^{t,x} = x$, and 
whose coefficients are assumed throughout this paper to
satisfy the following condition: 
\begin{itemize}
\item[(\hypertarget{A1}{$A_1$})]
  For every $t\in [0,T]$, 
 the functions 
$x\mapsto \mu(t,x)$ and $x\mapsto \sigma(t,x)$ are globally Lipschitz, i.e. 
\begin{equation} 
\nonumber 
	|\mu(t,x)-\mu(t,y)| \leq C|x-y|
	\quad
	\mbox{and}
	\quad 
	|\sigma(t,x)-\sigma(t,y)| \leq C|x-y|,
        \quad
        x,y\in \real, 
\end{equation}
\end{itemize}
In particular,
$x\mapsto \mu(t,x)$ and $x\mapsto \sigma(t,x)$ satisfy the linear growth conditions 
 \begin{equation} 
\nonumber 
 	|\mu(t,x)|  \leq C(1+|x|) \quad \mbox{and} \quad |\sigma(t,x)|  \leq C(1+|x|),
        \quad x\in \real, 
 \end{equation}
 for some positive constant $C > 0$. 
 The associated backward SDE is defined by
\begin{equation}
	\label{bw}
		Y_s^{t,x} = \phi\big(X_T^{t,x}\big) + \int_s^Tg\big( \tau ,X_\tau^{t,x}, Y_\tau^{t,x}, Z_\tau^{t,x} \big)d\tau
	- \int_s^TZ_\tau^{t,x} dB_\tau, \quad 0\leq t \leq s \leq T,
\end{equation}
with terminal condition $Y_T^{t,x}=\phi\big(X_T^{t,x}\big) \in L^2 (\Omega, \cF_T)$, 
where 
the generator $g(\cdot ,x,y, z)$ of \eqref{bw}
        is an $(\cF_t)_{t\in [0,T]}$-adapted process
        in $L^2\left(\Omega \times [0,T]\right)$
        for all $x,y,z\in \mathbb{R}$,
        which satisfies the following
 Conditions (\hyperlink{A2}{$A_2$})-(\hyperlink{A3}{$A_3$}). 
\begin{itemize}
\item[(\hypertarget{A2}{$A_2$})] The function $g(t,x,y,z)$ is uniformly Lipschitz
  in $(x,y, z)$, i.e., there exists $C > 0$ such that 
	\begin{equation} 
\nonumber 
		\left|g (t,x_2,y_2,z_2)-g(t,x_1,y_1,z_1)\right| \leq C\left(\left|x_2-x_1\right|+\left|y_2-y_1\right|+\left|z_2-z_1\right|\right) 
	\end{equation}
        a.s.,
        $x_1,x_2,y_1, y_2, z_1,z_2 \in \mathbb{R}$,
        $t\in [0,T]$, 
      \item[(\hypertarget{A3}{$A_3$})]
        We have $g(\cdot ,x,0,0) = 0$ a.s. for all $x \in \mathbb{R} $.
\end{itemize}
\noindent
By Theorem~2.1 of \cite{el1997backward},
see
Proposition~2.2 of \cite{pardoux1990adapted}, 
under
(\hyperlink{A1}{$A_1$})-(\hyperlink{A3}{$A_2$})  
there exists a unique pair $\big( Y_s^{t,x},Z_s^{t,x} \big)_{s\in [t,T]}$ 
of adapted processes in $L^2\left(\Omega \times [0,T]\right)$ 
that solves the BSDE \eqref{bw}.

\medskip

\noindent
 In the sequel we will also consider the following condition: 
\begin{itemize}
\item[(\hypertarget{A4}{$A_4$})]
 The function $\phi$ is continuous on $\real$ 
 and has the polynomial growth 
\begin{equation} 
  	\label{growth.phi}
 	|\phi(x)| \leq C(1+|x|^p), \qquad x\in \real,
        \mbox{~for some } p\geq1 \mbox{~and~} C>0.
\end{equation}
 \end{itemize}
 \subsubsection*{$g$-evaluation and $g$-expectation} 
 Next, we state the definition of the $g$-evaluation.
\begin{definition}
 Given $\xi \in L^2 (\Omega, \cF_T)$ and the backward SDE 
 \begin{align}
 \left\{
 \begin{array}{l}
 	dY^{0,x}_t =-g\big(t,X^{0,x}_t,Y^{0,x}_t,Z^{0,x}_t\big)dt + Z^{0,x}_t dB_t, \quad 0 \leq t\leq T,
 	\\
 	\\
 	Y^{0,x}_T = \xi, 
 \end{array}
 \right.
 \label{def_bsdes}
 \end{align}

 we respectively call
	\begin{equation}
          \nonumber
	  \cE_g[\xi] := Y^{0,x}_0
          \quad
          \mbox{and}
          \quad
          \cE_g[\xi \mid \cF_t] := Y^{0,x}_t
	\end{equation}
	the $g$-evaluation and the
        $\cF_t$-conditional $g$-evaluation of $\xi$,
             $t\in [0,T]$.
\end{definition}
Under (\hyperlink{A3}{$A_3$}),
one can show in addition that the map $\xi \mapsto \cE_g[\xi] $ preserves all properties of the classical expectation $\E$, except for linearity
and the property $\cE_g[c] = c$ for constant $c\in \real$,
see Relation~(34) and Theorem~3.4 in \cite{peng2004nonlinear}. 

\medskip

\noindent
 In the sequel we make the (stronger than (\hyperlink{A3}{$A_3$})) assumption 
\begin{itemize}
\item[(\hypertarget{A3'}{$A'_3$})] $g(\cdot ,x,y,0) = 0$ a.s.
  for all $x,y \in \mathbb{R}$, 
\end{itemize}
under which the $g$-evaluation $\cE_g$ becomes
the $g$-expectation, which satisfies
the property $\cE_g[c] = c$ for constant $c\in \real$,
see Relation~(36.2) and Lemma~36.3 in \cite{peng1997backward}.
We note that the results of Sections~\ref{s3}, \ref{s3.1} and
\ref{s5} remain valid for $g$-expectations if we assume
(\hyperlink{A3'}{$A'_3$}) instead of (\hyperlink{A3}{$A_3$}).
 
\medskip

\begin{remark}
  \begin{itemize}
  \item 
When $g(t,x,y,z)$ is convex in $ (y,z) \in \real^2$ for all
$ (t,x) \in [0,T]\times \real $,
which is the case in Theorems~\ref{pp3}-\ref{pp4}, 
Corollaries~\ref{pp5}-\ref{pp6}, 
and Theorems~\ref{pp1}-\ref{pp2}, we have the representation 
\begin{equation}
  \label{ddd}
  \cE_g [\xi] = \sup \limits_{\mathbb{Q} \in \mathcal{P}_g}
\left( \E_\mathbb{Q}[\xi] -F_g \left( \frac{d\mathbb{Q}}{d\mathbb{P}}\right)  \right),
\end{equation}
  	where $ F_g : L^2(\cF_T) \mapsto \real \cup \{+\infty\}$ is the convex
        functional defined by 
	$$ F_g (X) := \sup \limits_{\xi  \in L^2(\Omega,\cF_T)}
        ( \E [\xi X ] -\cE_g [\xi] ), 
        \quad X \in L^2(\Omega, \cF_T), 
        $$
	and $ \mathcal{P}_g$
        is the non-empty convex set
        of prior probability measures
        representing model uncertainty and defined by
	\begin{equation}
          \label{pgi}
          \mathcal{P}_g :=
        \left\{
        \mathbb{Q} \in \mathcal{M} \ \! : \ \! \frac{d\mathbb{Q}}{d\mathbb{P}} \in L^2(\cF_T) \mbox{ and } F_g \left(\frac{d\mathbb{Q}}{d\mathbb{P}}\right) < \infty \right\}, 
        \end{equation}
 where $\mathcal{M} $ is the set of probability measures
 on $ (\Omega, \cF_T) $ which are absolutely continuous with respect to $\mathbb{P}$,
 see Corollary~12 of \cite{gianin2006risk}. 
\item 
\textcolor{newcolor}{
 If $g (t,x,y,z)$ is both convex and sublinear in $ (y,z) \in \real^2$ for all $ (t,x) \in [0,T]\times \real $, which is the case in Examples~\ref{ex6.2}-\ref{ex6.4}, then \eqref{ddd}
 becomes 
 \begin{equation}
   \label{sas} 
 \cE_g [\xi] = \sup \limits_{\mathbb{Q} \in \mathcal{P}_g } \E_{\mathbb{Q}}[\xi], 
\end{equation} 
 see Corollary~12 in \cite{gianin2006risk}  and also 
 \cite{chen2000general}, \cite{chen2003jensen2}.
}
\item
  \textcolor{newcolor}{
      If $g(t,x,y,z)=\alpha(t,x) y +\beta(t,x)z$
 is linear in  $(y,z) \in \real^2$
      for all $ (t,x) \in [0,T]\times \real $,  
 where $\alpha (t,x)$ and $\beta (t,x)$
      are bounded functions, see Example~\ref{ex6.1},
      then the $g$-expectation
      $\cE_g [ \ \! \cdot \ \! ]$
      satisfies 
      \begin{equation} 
\nonumber 
 	\cE_g\big[\xi~\big| \cF_t\big]=\E_\mathbb{Q}\left[
          \xi \exp \left( \int_t^T\alpha\big(s,X_s\big) ds \right)
          \ \! \Big| \ \! \cF_t\right], 
 \end{equation} 
      where
      $\mathbb{Q}$
      is the probability measure defined as 
\begin{align}
 	 \frac{d\mathbb{Q}}{d\mathbb{P}} := 
 	\exp\left(\int_0^T \beta (s,X_s) dB_s -\frac{1}{2}\int_0^T \beta^2 (s,X_s)ds\right), 
 \end{align}
as follows by applying the It\^{o} formula 
to $Y_t\exp\big(\int_0^t\alpha(s,X_s)ds\big)$
using the Brownian motion $\widetilde{B}_t:=B_t+\int_0^t\beta(s,X_s)ds$ under
 $\mathbb{Q}$.
}
     \item
  \textcolor{newcolor}{
    If $g(t,x,y,z) = \alpha_t |z| + \beta_t z$,
    $t\in [0,T]$, $z \in \real$, where
    $(\alpha_t)_{t\in [0,T]}$ and $(\beta_t )_{t\in [0,T]}$ are
    time-continuous processes then 
    $\cE_g[\xi]$ coincides with Choquet's expectation
\begin{align}
 \nonumber 
 \cE_g[\xi] = \E_{\mu_g}[\xi]
  :=\int_{-\infty}^{0}(\mu_g (\xi>x)-1)dx + \int_0^{\infty}\mu_g (\xi>x)dx
\end{align}
for $\xi$ of the form $\xi=y+zB_T$, where $\mu_g[A]:=\cE_g[\mathbf{1}_A]$,
 $A \in \cF$, is the corresponding capacity,
    see \cite{chen-davison}.
   } 

  \medskip

  \textcolor{newcolor}{
    Moreover, $\cE_g [ \ \! \cdot \ \! ]$
    coincides with the linear expectation
    $\cE_{\mu_g}[ \ \! \cdot \ \! ]$
    if and only if $\alpha_t =0$ $a.s.$, $t\in [0,T]$, 
    see Theorem~1 in \cite{chen-davison},
    and in case $\beta_t=0$, $t \in [0,T]$, we have
    \vspace{-0.5cm} 
  \begin{subequations}
          \begin{empheq}[left={\cE}_g [\xi{]}{=}\empheqlbrace]{align}
\label{g_alpha-a} 
 	 		\sup\limits_{\mathbb{Q} \in \mathcal{P}_g } \E_\mathbb{Q}[\xi], & \mbox{ if } \alpha_t >0, \ a.s.,  \ t\in [0,T],
        		\\
        \label{g_alpha-b} 
 		\inf\limits_{\mathbb{Q} \in \mathcal{P}_g }\E_\mathbb{Q}[\xi], & \mbox{ if } \alpha_t <0, \ a.s., \ t\in [0,T],
        	\end{empheq}
  \end{subequations}
  where $\mathcal{P}_g$ is the set of probability measures
  $\mathbb{Q} \in \mathcal{M}$ such that 
$$
 \frac{d\mathbb{Q}}{d\mathbb{P}}:= \exp \left( -\int_0^Tv_tdB_t
 - \frac{1}{2} \int_0^Tv_t^2dt\right), 
$$
 where $(v_t)_{t\in [0,T]}$ is $(\cF_t)_{t\in [0,T]}$-adapted and
  $|v_t| \leq |\alpha_t|$, $t\in [0,T]$, 
see Theorem~2.2 in \cite{chen2002ambiguity} and also Example 1 in \cite{chen-davison}.
}  
\end{itemize}
\end{remark}
\noindent
In the sequel, we state ordering results for the $g$-evaluation 
$\cE_g[ \ \! \cdot \ \! ]$, and more generally for the 
conditional $g$-evaluation and the conditional $g$-expectation
$\cE_g[ \ \! \cdot \mid \cF_t]$, $t\in [0,T]$,
under (\hyperlink{A3'}{$A'_3$}) instead of (\hyperlink{A3}{$A_3$}).
\subsubsection*{$g$-stochastic orderings}
Stochastic orderings with respect to capacity have been studied in \cite{grigorova2014stochastic,grigorova2014dominance} by using Choquet's expectation and uncertainty orders have been constructed in \cite{tian2016uncertainty}
on the sublinear $G$-expectation space. 
Here, we extend their approaches to the comparison of
random variables $X^{(1)}$, $X^{(2)}$ in the 
 settings of Peng's $g$-expectations 
and $g$-evaluations, which are not sublinear in general,
via the condition
\begin{equation}
          \label{co2} 
	          \cE_{g_1}\big[ \phi\big(X^{(1)}\big)\big] \leq \cE_{g_2}\big[ \phi\big(X^{(2)}\big)\big],
		\end{equation}
                in nonlinear expectations $\cE_{g_1}$, $\cE_{g_2}$,  
        for all $\phi(x)$
        in a certain class of functions having polynomial growth.
        In general, different portfolios or hedging strategies
may corresponding to different generators $g_1$, $g_2$
as can be seen in Examples~\ref{ex6.2} and \ref{ex6.4}.
\begin{definition}
    Let $g_1,g_2$ satisfy (\hyperlink{A2}{$A_2$})-(\hyperlink{A3}{$A_3$}).
    For any $X^{(1)}, X^{(2)} \in L^2 (\Omega, \cF_T )$, we say that 
	\begin{enumerate}[1)]
	\item $X^{(1)}$ is dominated by $X^{(2)}$
          in the monotonic $g_1,g_2$-ordering, i.e. $X^{(1)} \le_{g_1,g_2}^{\rm mon} X^{(2)}$,
          if \eqref{co2} 
          holds for all non-decreasing functions $\phi(x)$
                  satisfying \eqref{growth.phi}. 
	      \item $X^{(1)}$ is dominated by $X^{(2)}$ in the convex $g_1,g_2$-ordering,
                i.e. $X^{(1)} \le_{g_1,g_2}^{\rm conv} X^{(2)}$, 
                if \eqref{co2} holds for all convex functions
                $\phi(x)$ satisfying \eqref{growth.phi}. 
		              \item $X^{(1)}$ is dominated by $X^{(2)}$ in the increasing convex $g_1,g_2$-ordering, i.e. $X^{(1)} \le_{g_1,g_2}^{\rm iconv} X^{(2)}$, if \eqref{co2} holds
                                for all non-decreasing convex functions
                                $\phi(x)$ satisfying \eqref{growth.phi}. 
	\end{enumerate}	
\end{definition}
We note that 
\begin{enumerate}[(i)]
\item $X^{(1)} \le_{g_1,g_2}^{\rm mon} X^{(2)} \Longrightarrow X^{(1)} \le_{g_1,g_2}^{\rm iconv} X^{(2)}$, and 
		\item $X^{(1)} \le_{g_1,g_2}^{\rm conv} X^{(2)} \Longrightarrow X^{(1)} \le_{g_1,g_2}^{\rm iconv} X^{(2)}$, 
	\end{enumerate}
and we simply write 
$\le_g^{\rm mon}$, $\le_{g}^{\rm conv}$, $\le_g^{\rm iconv}$
if $g_1=g_2:=g$.
\textcolor{newcolor}{Replacing convex functions with concave functions
  yields the corresponding notions
  of concave and increasing concave $g_1,g_2$-orderings
  denoted by $X^{(1)} \le_{g_1,g_2}^{\rm conc} X^{(2)}$ and $X^{(1)} \le_{g_1,g_2}^{\rm iconc} X^{(2)}$, respectively,
  which are characterized in the next proposition. 
  \begin{prop}
    \label{fdsf} 
	Given two random variables $X^{(i)}\in L^2(\Omega, \cF_T), \, i=1,2$ and the generator $ g^{(-1)}(t,x,y,z) :=-g(t,x,-y,-z),$ for $ (t,x,y,z) \in [0,T] \times \real^3 $, we have that
	\begin{enumerate}[(i)]
		\item $X^{(1)} \le_g^{\rm conc} X^{(2)}$ if and only if $-X^{(2)} \le_{g^{(-1)}}^{\rm conv} -X^{(1)}$.
		\item $X^{(1)} \le_g^{\rm iconc} X^{(2)}$ if and only if $-X^{(2)} \le_{g^{(-1)}}^{\rm iconv} -X^{(1)}$.
	\end{enumerate} 
\end{prop}
  Proposition~\ref{fdsf} is
  a direct consequence of the following lemma.
}
\textcolor{newcolor}{
  \begin{lemma}
  \label{lem}
Letting $g^{(a)}(t,x,y,z) := ag(t,x,y / a, z / a )$,
	$x,y,z\in \real$, $t\in [0,T]$, we have  
	$$
	\cE_{{g}^{(a)}}[ a\xi]
	= a\cE_{{g}^{(1)}} [ \xi ], 
        \qquad a \ne 0, 
        $$  
\end{lemma}
}
\begin{Proof}
              Given $(X_t)_{t\in[0,T]}$ an
     $(\cF_t)_{t\in [0,T]}$-adapted process, 
	let $(Y_t,Z_t)_{t\in [0,T]}$ denote the solution of 
        \eqref{def_bsdes} 
	with terminal condition $Y_T=\xi
	\in L^2(\Omega, \cF_T)$, and let 
	$(\widebar{Y}_t,\widebar{Z}_t)_{t\in [0,T]}$
	denote the solution of the backward SDE
	\begin{equation}
		\nonumber
		\widebar{Y}_t = a\xi + \int_t^T {g}^{(a)} \left(s, X_s, \widebar{Y}_s, \widebar{Z}_s\right) ds - \int_t^T \widebar{Z}_s dB_s,
		\qquad t\in [0,T],
	\end{equation}
	with generator
	$g^{(a)} (t,x,y,z) $ and terminal condition $\widebar{Y}_T = a\xi $,
	i.e.
	$$ 
	\frac{\widebar{Y}_t }{a} = \xi + \int_t^T g \left(s, X_s, \frac{\widebar{Y}_s}{a}, \frac{\widebar{Z}_s}{a}\right) ds
	- \frac{1}{a} \int_t^T \widebar{Z}_s dB_s,
	\qquad t\in [0,T].
	$$ 
	The uniqueness of the solution $(Y_t, Z_t)$ of the backward SDEs
	\begin{equation}
		\nonumber 
		Y_t = \xi + \int_t^T g\left(s, X_s, Y_s, Z_s\right) ds - \int_t^TZ_s dB_s, 
		\qquad t\in [0,T], 
	\end{equation}
	yields $\widebar{Z}_t / a  = Z_t$ and
        $\widebar{Y}_t / a = Y_t$, $t\in [0,T]$, and 
        we conclude by taking $t=0$.
	          \end{Proof}
\noindent 
We also note that the monotonic $g$-ordering admits the following characterization
in the case of sublinear generator functions.~ 
\begin{prop} Assume that $ g_i(t,x,y,z)$ 
  is sublinear in $ (y,z) \in \real^2$ for all $ (t,x) \in [0,T]\times \real $, $i=1,2$, 
  and that $ g_1(t,x,y,z) \leq g_2(t,x,y,z) $.
    Then $ X^{(1)} \le_{g_1,g_2}^{mon}  X^{(2)} $ is equivalent to
  \begin{equation}
    \label{sh} 
  \inf \limits_{Q \in \mathcal{P}_{g_1}}Q(X^{(1)} \leq c) \ge \inf \limits_{Q \in \mathcal{P}_{g_2}}Q(X^{(2)} \leq c), \qquad
  c \in \real ,
\end{equation} 
  where
 $\mathcal{P}_{g_1}$ and $\mathcal{P}_{g_2}$ are defined in \eqref{pgi}. 
\end{prop}
	\begin{Proof} 
	  {$(i) \Rightarrow (ii)$:} We apply \eqref{sas}
          to the non-decreasing function
          $ \phi(x):= \mathbf{1}_{\{x > c \}}$
for $c \in \real$,
after noting that $ \mathcal{P}_{g_1} \subseteq \mathcal{P}_{g_2}$
since $ g_1 \leq g_2 $ by Remark~13 in \cite{gianin2006risk}. 
\\
{$(ii) \Rightarrow (i)$:}
By relation \eqref{sh}, for any $ Q \in \mathcal{P}_{g_1}$ and non-decreasing functions
$\phi $ we have 
$$  \E_Q\big[\phi\big(X^{(1)}\big)\big] \leq \E_Q\big[\phi\big(X^{(2)}\big)\big] \leq \sup \limits_{Q \in \mathcal{P}_{g_2}} \E_Q\big[\phi\big(X^{(2)}\big)\big], 
$$ 
hence by \eqref{sas} we find 
$$
\cE_{g_1}\big[\phi\big(X^{(1)}\big)\big]
=
\sup \limits_{Q \in \mathcal{P}_{g_1}} \E_Q\big[\phi\big(X^{(1)}\big)\big]
\leq \sup \limits_{Q \in \mathcal{P}_{g_2}} \E_Q\big[\phi\big(X^{(2)}\big)\big]
=
\cE_{g_2}\big[\phi\big(X^{(2)}\big)\big].
$$ 
\end{Proof}
\subsubsection*{Associated PDE}
\noindent
Throughout the remaining of this paper
we assume that $g(t,x,y,z)$ is a deterministic function,
in addition to 
(\hyperlink{A1}{$A_1$})-(\hyperlink{A3}{$A_3$}).
 The function $ u(t,x):=Y_t^{t,x} $
 can be shown to be 
 a viscosity solution of the backward PDE
   \begin{equation}
	\label{du12} 
  	\frac{\partial u}{\partial t}(t,x) +\mu(t,x) \frac{\partial u}{\partial x}(t,x) 
	+ \frac{1}{2}\sigma^2(t,x) \frac{\partial^2 u}{\partial x^2}(t,x)+ g\Big(t,x,u(t,x),\sigma(t,x) \frac{\partial u}{\partial x}(t,x) \Big)=0,
\end{equation}
   $x\in \real$, $t\in [0,T]$, with a terminal condition $u(T,x)=\phi(x)$
   satisfying (\hyperlink{A4}{$A_4$}), see Theorem~2.2 in \cite{pardoux1998backward},
   Theorem~4.3 of \cite{pardouxpeng} and Theorem~4.2 of \cite{el1997backward}. 

\medskip

In the sequel, we let ${\cal C}^{p,q}([0,T]\times\real)$ 
denote the space of functions $f(t,x)$ which are
$p$ times continuously differentiable in $t\in [0,T]$,
$p \geq 1$, 
and $q$ times differentiable in $x\in \real$, 
$q\geq 1$. 
We also let ${\cal C}_b^k(\real^n)$ 
denote the space of continuously differentiable functions 
   whose partial derivatives of orders one to $k$ are uniformly bounded
   on $ \real^n $. 
In Theorem~\ref{classical} below we state an existence result for
    classical solutions under stronger smoothness assumptions on
    BSDE coefficients, see Theorem~3.2 of \cite{pardouxpeng},
    Theorem~8.1 in \S~V.8 page~495, and
    Theorem~7.1 in \S~VII.7 page~596 
of \cite{ladyzenskaja}. 
   \begin{theorem}\label{classical}
   Assume (\hyperlink{A3}{$A_3$}) 
   and in addition that $ \mu(t,\cdot )$,
   $\sigma(t,\cdot )$, $\phi \in {\cal C}_b^3(\real)$,
   and that $ g(t,\cdot ,\cdot ,\cdot )\in {\cal C}_b^3 (\real^3) $ 
   for any $ t\in [0,T] $. 
   Then the function $u(t,x):=Y_t^{t,x}$ is a classical solution in
   ${\cal C}^{1,2}([0,T]\times\real)$ of the backward PDE
   \eqref{du12} with 
  terminal condition $u(T,\cdot )=\phi$. 
\end{theorem}
 Under the conditions of Theorem~\ref{classical},
 by Proposition~4.3 of \cite{el1997backward}
        the solution $\big(Y^{t,x}_s,Z^{t,x}_s\big)_{s\in [t,T]}$
        of \eqref{abs}
        satisfies 
        $Y^{t,x}_s = u (s,X^{t,x}_s) $ 
        and
        $\displaystyle Z^{t,x}_s = \sigma (s,X^{t,x}_s) \frac{\partial u}{\partial x} (s,X^{t,x}_s)$,
$0\leq t \leq s \leq T$. 
         In addition, by Theorem~2.2 in \cite{ma1999forward}
 or Proposition~3.3 in \cite{maprotteryong}
 we have the following result. 

  \begin{theorem}\label{classical2}
   Under the assumptions of Theorem~\ref{classical},
   suppose additionally that $\sigma (t,x)$ is
   bounded above and below by strictly positive constants.
Then the first derivative in $t\in [0,T]$ 
 and the first and second derivatives in $x\in \real$ of $u(t,x)$
 are bounded in $(t,x)\in [0,T]\times \real$. 
  \end{theorem}
  As in \cite{DMP96}, we denote by
  ${\cal C}^{1+\eta/2,2+\eta}([0,T]\times \real)$, $\eta \in (0,1)$,
 the space of functions $f(t,x)$ which are
differentiable in $t\in [0,T]$ 
and twice differentiable in $x\in \real$
with
$\frac{\partial f}{\partial t}(t,x)$ and 
$\frac{\partial^2 f}{\partial x^2}(t,x)$ being
respectively $\eta/2$-H\"older continuous and 
$\eta$-H\"older continuous in $(t,x)\in [0,T]\times \real$,
and define the space ${\cal C}^{k+\eta }( \real)$ analogously for $k\geq 1$.
By Theorem~2.3 in \cite{DMP96}, see also page 236 of
\cite{ma1999forward}, we have the following result.
 
\begin{theorem}\label{classical3}
   In addition to the assumptions of Theorem~\ref{classical2},
   suppose that  for some $\eta \in (0,1)$
   the functions $\mu (\cdot, \cdot)$, $\sigma (\cdot, \cdot)$ and $g(\cdot ,\cdot ,y,z)$ 
   are in ${\cal C}^{1+\eta /2,2+\eta }([0,T]\times \real)$
   for all $y,z\in \real$, and
    that $\phi \in{\cal C}^{4+\eta }( \real)$.
   Then the function $u(t,x)$ is a classical solution in
   ${\cal C}^{2+\eta / 2,4+\eta }([0,T]\times\real)$ of the backward PDE
   \eqref{du12} with 
  terminal condition $u(T,\cdot )=\phi$. 
 \end{theorem} 

\section{Ordering with convex drifts} 
                        \label{s3}
                       Consider the forward SDEs
        \begin{subequations}
        	\begin{empheq}[left=\empheqlbrace]{align}
        		\label{d3} 
        	&	dX_t^{(1)}=\mu_1\big(t,X_t^{(1)}\big)dt + \sigma_1\big(t,X_t^{(1)}\big)dB_t, 
        		\\
        		\nonumber 
        		\\ 
    		\label{d4} 
        	&	dX_t^{(2)}=\mu_2\big(t,X_t^{(2)}\big)dt + \sigma_2\big(t,X_t^{(2)}\big)dB_t, 
        	\end{empheq}
        \end{subequations}
    and  the associated BSDEs 
    \begin{equation}
      \label{abs} 
    \left\{
    \begin{array}{ll}
    \nonumber
    dY_t^{(1)} = - g_1\big(t,X_t^{(1)},Y_t^{(1)},Z_t^{(1)}\big) dt + Z_t^{(1)}dB_t, & Y_T^{(1)} =\phi\big( X_T^{(1)}\big),
    \\
    \\
    \nonumber
    dY_t^{(2)} = -g_2\big(t,X_t^{(2)},Y_t^{(2)},Z_t^{(2)}\big) dt + Z_t^{(2)} dB_t, & Y_T^{(2)} =\phi\big( X_T^{(2)}\big),  
    \end{array}
    \right.
\end{equation} 
    and let
    \begin{equation}
      \label{fi} 
    f_i ( t,x,y,z): = 
      	z\mu_i (t,x)
  	+
  	g_i (t,x,y,z \sigma_i ),
        \quad
        t\in [0,T],
        \  
        x,y,z\in \real,
\         i = 1,2.
\end{equation} 
        In all following propositions, the 
convexity of
$$(x,y) \mapsto f_i(t,x,y,z),
\quad
\mbox{resp.}
\quad 
(y,z) \mapsto f_i(t,x,y,z)
$$
on $\real^2$ is understood to hold for all $(t,z)\in [0,T]\times \real$,
resp.
for all $(t,x)\in [0,T]\times \real$.
The next result is a consequence of the 
   comparison Theorem~2.4 in Appendix~C of \cite{peng}. 
   We note that
   Condition~(\hyperlink{B1}{$B_1$}) can be shown to be
   necessary for convex ordering by taking $\phi (x) = x$
   as in Theorem~3.2 of \cite{briand2000converse}. 
   \begin{theorem}
  	\label{pp3}
        {\em (Convex order)}.
  	Assume that $X_0^{(1)}= X_0^{(2)}$, and  
                $$ 0 < \sigma_1(t,x) \leq \sigma_2(t,x), \qquad t\in [0,T], \quad x \in \real, $$ together with the conditions  
\begin{itemize} 
\item[{\em (\hypertarget{B1}$B_1$)}]
$	\displaystyle  
  	f_1(t,x,y,z)
  	\leq
  	f_2(t,x,y,z)
$, $t\in [0,T]$, $x,y,z \in \real$,
\item[{\em (\hypertarget{B2}$B_2$)}]
        $\displaystyle (x,y) \mapsto f_i(t,x,y,z)$ and $(y,z) \mapsto f_i(t,x,y,z)$
  are convex on $\real^2$ for $i=1,2$.
  	\end{itemize} 
  	Then we have $X_T^{(1)}\le_{g_1,g_2}^{\rm conv} X_T^{(2)}$,  i.e.,
        \begin{equation}
          \label{ceg0} 
  	\cE_{g_1}\big[\phi\big(X_T^{(1)}\big)\big] \leq \cE_{g_2}\big[\phi\big(X_T^{(2)}\big)\big],
\end{equation} 
  	for all convex functions $\phi(x)$ satisfying \eqref{growth.phi}. 
  \end{theorem}
   \begin{Proof} 
 We start by assuming that
 the function $\phi$ and the coefficients $\mu_i(t,\cdot)$, $\sigma_i(t,\cdot)$
 and $g_i(t,\cdot,\cdot,\cdot)$ are 
 ${\cal C}_b^3$ functions for all $t\in [0,T]$.
 By Theorem~\ref{classical},
the functions $u_1(t,x) := Y_t^{(1),t,x} $ and $u_2(t,x) := Y_t^{(2),t,x}$
are 
solutions of the backward PDEs
	\eqref{du1-1} which are continuous in $t $ and $x$.  
         Letting
        \begin{equation}
          \label{dsdsas} 
	h_i(t,x,y,z,w) := f_i(t,x,y,z)
        + \frac{w}{2}\sigma^2_i(t,x),
	\qquad
	i=1,2, 
\end{equation} 
we rewrite \eqref{du12} 
	as 
	\begin{equation}
        \label{du1-1} 
	\frac{\partial u_i}{\partial \tau}(\tau,x)
        =
        h_i\left(\tau, x, u_i(\tau,x),\frac{\partial u_i}{\partial x}(\tau,x), \frac{\partial^2 u_i}{\partial x^2}(\tau,x)\right)
        \text{ with } u_i(0,x) = \phi(x),
	\quad
	i=1,2,
\end{equation} 
        by setting $\tau:=T-t $. 
               We also assume that there exists constants $c,C'>0$ such that
        \begin{equation}
          \label{sigma}
          0 < c \leq \sigma_1(t,x) \leq \sigma_2(t,x) \leq C', \qquad t\in [0,T], \quad x \in \real.
        \end{equation} 
        In this case, by Theorem~\ref{classical2}
        the second derivative
        $\left| \displaystyle \frac{\partial^2 u_i}{\partial x^2}(t,x)\right|$
        is
        bounded by $C''>0$. 
	In addition, under (\hyperlink{B1}{$B_2$}),
        both solutions $u_1(t,x) $ and $u_2(t,x) $ of \eqref{du1-1}
        are convex functions of $x$ by Theorem~\ref{convex} below, 
	hence we have
        $\displaystyle \frac{\partial^2 u_i}{\partial x^2}(\tau,x)\geq 0$,
$\tau \in [0,T]$, $x\in \real$.
Therefore, in \eqref{du1-1} we 
        can replace $h_i(t,x,y,z,w)$ in \eqref{dsdsas} with 
	\begin{equation} 
        \label{dsdsas2}
 	\tilde{h}_i(t,x,y,z,w) := f_i (t,x,y,z)
        + \frac{(\min ( w , C''))^+}{2}\sigma^2_i(t,x),
	\qquad
	i=1,2,
\end{equation} 
	where $w^+=\max (w,0)$,
        and rewrite the backward PDEs \eqref{du1-1} as 
	$$ 
	\frac{\partial u_i}{\partial \tau}(\tau,x)
        =
        \tilde{h}_i\left(\tau, x, u_i(\tau,x),\frac{\partial u_i}{\partial x}(\tau,x), \frac{\partial^2 u_i}{\partial x^2}(\tau,x)\right)
        \text{ with } u_i(0,x) = \phi(x),
	\quad
	i=1,2. 
	$$
        Next, for all $\tau \in [0,T]$ and $x_1,x_2,y,z\in \real$ we have 
        	\begin{eqnarray*}
		\lefteqn{
  \! \! \! \! \! \! \! \! \! \! \! \! \! \! \! 
  |              f_i(\tau, x_2,y,z)-f_i(\tau, x_1,y,z) |
  		}
		\\
		& \leq & 
  |z|
  |\mu_i(\tau, x_2)-\mu_i(\tau, x_1)
			| 
+ | g_i(\tau, x_2,y, z\sigma_i(\tau, x_2))-g_i(\tau, x_1,y, z\sigma_i(\tau, x_1))|
		\\
		&\leq & C |z| |x_2-x_1|  + C(|x_2-x_1|+|z|
                \left|\sigma_i(\tau, x_2)-\sigma_i(\tau, x_1)\right|)
                \\
		&\leq & C |z| |x_2-x_1|  + C ( 1 + |z| ) |x_2-x_1|, \qquad
		i=1,2, 
	\end{eqnarray*}
 hence 
        	\begin{eqnarray*}
		\lefteqn{
  \! \! \! 
  |              \tilde{h}_i(\tau, x_2,y,z,w)-\tilde{h}_i(\tau, x_1,y,z,w) | 
		}
		\\
		& \leq &
                 | f_i(\tau, x_2,y, z)-f_i(\tau, x_1,y, z)|
                 +  \frac{(\min ( w , C''))^+}{2}
                 |\sigma_i^2(\tau, x_2)-\sigma_i^2(\tau, x_1)|
                		\\
		                &\leq &
                                C |z| |x_2-x_1|  + C ( 1 + |z| ) |x_2-x_1|
                                +
                                \frac{C''}{2}
                                \left|\sigma_i(\tau, x_2)-\sigma_i(\tau, x_1)\right|
                 ( \sigma_i(\tau, x_1)+\sigma_i(\tau, x_2) ) 
\\
&\leq & C |z| |x_2-x_1| +  C(1+|z|) | x_2- x_1| + C C'C'' |x_2-x_1| 
           \\
	   & \leq & ( C + CC'C'') \left( 1 
           + |x_1| + |x_2| +|y| \right)
		( 1 + |z| ) |x_2-x_1 |, \quad
		i=1,2, 
	\end{eqnarray*}
                	        which shows that Condition~(G) of Theorem~2.4 in Appendix~C of \cite{peng}
        is satisfied with $\omega(x)=\bar{\omega}(x):=Cx$.
        In addition, by the conditions \eqref{sigma} and
                (\hyperlink{B1}{$B_1$}) we have
	\begin{eqnarray*}
          \lefteqn{
            \! \! \! \! \! \! \! \! \! \! \! \! \! 
            	    \tilde{h}_2(\tau,x,y,z,w)- \tilde{h}_1(\tau,x,y,z,w)
}
\\
           & = & 
               f_2(\tau,x,y,z )
- f_1(\tau,x,y,z ) 
+ \frac{(\min ( w , C''))^+}{2}\big(
	  \sigma^2_2(\tau,x)-\sigma^2_1(\tau,x)\big) 
		\\
		& \ge & 0,
                \qquad
x,y,z,w \in \real, \quad \tau \in [0,T]. 
	\end{eqnarray*}
	Besides, we have
	$\tilde{h}_2(\tau,x,y,z,w_1) \leq \tilde{h}_2(\tau,x,y,z,w_2)$ 
        when $w_1 \leq w_2 $, and 
\begin{eqnarray*}
  \lefteqn{
    \! \! \! \! \! \! \! \! \! \! \!
    \! \! \! \! \! \! \! \! \! \! \!
    |\tilde{h}_2(\tau,x,y_1,z,w_1)- \tilde{h}_2(\tau,x,y_2,z,w_2)| \leq \frac{1}{2}\sigma^2_2(\tau,x)\left|w_1-w_2\right|
  }
  \\
& & 
+\left|f_2(\tau,x,y_1,z)-f_2(\tau,x,y_2,z)\right|
\\
& \leq & C\left(|y_1-y_2| +|w_1-w_2| \right), 
\end{eqnarray*}
$ (\tau,x) \in [0,T]\times \real $,
$ (y_1,z,w_1)$, $(y_2,z,w_2) \in \real^3$,
hence $\tilde{h}_2(t,x,y,z,w) $ is Lipschitz in $y $ and $w$.
Therefore, by the comparison Theorem~2.4 in Appendix~C of \cite{peng} 
	it follows that $u_1(t, x) \leq u_2(t,x)$
	for all $t\in [0,T]$, from which we conclude to
        $$
        Y_0^{(1)} = u_1\big(0, X_0^{(1)}\big) \leq Y_0^{(2)} = u_2\big(0,X_0^{(1)}\big) =u_2\big(0,X_0^{(2)}\big),
 $$ 
	hence $$\cE_{g_1} \big[ \phi\big(X_T^{(1)} \big)\big]
	\leq  \cE_{g_2} \big[ \phi \big( X_T^{(2)} \big)\big],
        $$ 
	for all convex functions $\phi$ in ${\cal C}^3_b(\real )$.
                In order to extend \eqref{ceg0} to coefficients satisfying
                                (\hyperlink{A1}{$A_1$})-(\hyperlink{A4}{$A_4$}) 
without assuming the bound \eqref{sigma}, 
we apply the above argument to sequences
$(\mu_{n,i})_{n\geq 1}$,
$(\sigma_{n,i})_{n \geq 1}$, 
$(g_{n,i})_{n\geq 1}$,
$(\phi_n)_{n\geq 1}$ of ${\cal C}_b^3$ functions
as in Theorem~\ref{classical},
with 
        \begin{equation}
          \label{sigma2}
          0 < c_n \leq \sigma_{n,1} (t,x) \leq \sigma_{n,2} (t,x) \leq C_n, \qquad t\in [0,T], \quad x \in \real, \quad n \geq 1,  
        \end{equation} 
        for some constants $c_n,C_n>0$ 
         satisfying (\hyperlink{A1}{$A_1$})-(\hyperlink{A4}{$A_4$})
 with a same constant $C>0$ for all $n\geq 1$,
 and converging respectively pointwise to $\mu_i$, $\sigma_i$, $g_i$  
 and strongly to $\phi$, i.e. $\phi_n(x_n) \to \phi(x)$
 whenever $x_n \to x \in \real$,
 while preserving the convexity of 
the approximations $(\phi_n )_{n\geq 1}$
and $(f_{n,i} )_{n\geq 1}$ defined by \eqref{fi},
see \cite{azagra}, 
Lemma~1 of \cite{lepeltier-sanmartin},
and Problem~1.4.14 in \cite{zhangjianfeng}. 
The continuous dependence Proposition~\ref{dependence} then yields the convergence
of the corresponding sequences $(Y^{(i)}_{n,0})_{n\geq 1}$
of BSDE solutions, concluding the proof. 
   \end{Proof} 
By similar arguments, 
we derive the following
Theorem~\ref{pp4} 
        for the increasing convex ordering.
        The proof of Theorem~\ref{pp4}
        is first stated for ${\cal C}_b^3$ coefficients
        $\phi$,
        $\mu_i(t,x)$,
$\sigma_i(t,x)$
        and $g_i(t,x,y,z)$ under \eqref{sigma2},
        and then extended to coefficients satisfying
 (\hyperlink{A1}{$A_1$})-(\hyperlink{A4}{$A_4$}) 
by applying the continuous dependence Proposition~\ref{dependence}
as in the proof of Theorem~\ref{pp3}. 
       \begin{theorem}
	\label{pp4}
        {\em (Increasing convex order)}. 
        Assume that $X_0^{(1)}\leq X_0^{(2)}$ and 
	$$
	0 < \sigma_1(t,x) \leq \sigma_2(t,x), \qquad
        t \in [0,T],
	\quad x \in \real,
	$$
	together with the conditions  
        \begin{itemize}
          \item[{\em (\hypertarget{B1}$B'_1$)}]
$
	\displaystyle  
  	f_1(t,x,y,z)
  	\leq
  	f_2(t,x,y,z)
$, 
 $t\in [0,T]$, $x,y \in \real$, $z\in \real_+$,
      \item[{\em (\hypertarget{B2}$B'_2$)}]
      $   	\displaystyle
   (x,y) \mapsto f_i(t,x,y,z)$  and $(y,z) \mapsto f_i(t,x,y,z)$
   are both convex respectively on $\real^2$ and $\real\times \real_+$, 
    for $i=1,2$, $x,y \in \real$, $z \in \real_+$, $t\in [0,T]$, 
      \item[{\em (\hypertarget{B3}$B'_3$)}] 
$\   	\displaystyle
   x \mapsto g_i(t,x,y,z)$
    is non-decreasing on $\real$ for $i=1,2$,
  $y \in \real$, $z \in \real_+$, $t\in [0,T]$.
\end{itemize} 
	Then we have $X_T^{(1)}\le_{g_1,g_2}^{\rm iconv} X_T^{(2)}$, i.e.,
	$$
	\cE_{g_1}\big[\phi\big(X_T^{(1)}\big)\big] \leq \cE_{g_2}\big[\phi\big(X_T^{(2)}\big)\big],
	$$ 
	for all non-decreasing convex functions $\phi(x)$ satisfying \eqref{growth.phi}. 
\end{theorem}
       \begin{Proof} 
         Under (\hyperlink{B'3}{$B'_3$}), 
         when $ \phi(x) $ and $ g_i(t,x,y,z)$, $i=1,2$, are non-decreasing in $ x $, 
Proposition~\ref{nondecreasing} tells us that the PDE solutions $u_1(t,x) $
	and $u_2(t,x) $
	satisfy 
	$$
	\frac{\partial u_1}{\partial x} (t,x) \ge 0 \
	\quad
	\mbox{and}
	\quad
	\frac{\partial u_2}{\partial x} (t,x) \ge 0,
        \qquad
        t\in [0,T], 
	$$
	hence Conditions~(\hyperlink{B1}{$B_1$})-(\hyperlink{B2}{$B_2$}) 
        only need to hold for $z\geq 0$, 
                and the conclusion follows by repeating the
                arguments in the proof of Theorem~\ref{pp3}. 
\end{Proof}
       We note that in case $ \sigma_1(t,x)=\sigma_2(t,x)$
       the convexity of $ u_i(t,x)$,
       $i=1,2$, is no longer required in the proofs of Theorems~\ref{pp3}-\ref{pp4},
       and one can then remove Condition~(\hyperlink{B'2}{$B'_2$})
       to obtain a result for the monotonic order. 
\begin{corollary}
	\label{pp4.1} 
              {\em (Monotonic order with equal volatilities)}.  
	 Assume that $X_0^{(1)} \leq X_0^{(2)}$ and $$ 0 < \sigma(t,x): = \sigma_1(t,x)=\sigma_2(t,x), \quad t\in [0,T], \quad x\in \real,
	$$ 
	together with the conditions 
 \begin{itemize} 
\item[{\em (\hypertarget{B''1}$B''_1$)}] 
$	\displaystyle  
  	f_1(t,x,y,z)
  	\leq
  	f_2(t,x,y,z)
$, 
 $t\in [0,T]$, $x,y \in \real$,
 $z\in \real_+$,
\item[{\em (\hypertarget{B''2}$B''_2$)}] 
        $\displaystyle
  x \mapsto g_i(t,x,y,z)$ is non-decreasing on $\real$ for $i=1,2$,
  $y \in \real$, $z \in \real_+$, $t\in [0,T]$.
\end{itemize} 
	Then we have $ X_T^{(1)}\le_{g_1,g_2}^{\rm mon} X_T^{(2)} $,  i.e.,
	$$
	\cE_{g_1}\big[\phi\big(X_T^{(1)}\big)\big] \leq \cE_{g_2}\big[\phi\big(X_T^{(2)}\big)\big],
	$$ 
	for all non-decreasing functions $\phi(x)$ satisfying \eqref{growth.phi}. 
\end{corollary}
\begin{Proof} 
  When $ \sigma_1(t,x) =\sigma_2(t,x)$ we can repeat the proof
  of Theorem~\ref{pp3} by using $h_i$ in \eqref{dsdsas}, 
  without defining $\tilde{h}_i$ in \eqref{dsdsas2}
  and without assuming (\hyperlink{B2}{$B_2$}),
  and then follow the proof argument of Theorem~\ref{pp4}
  without requiring the convexity of $ u_i(t,x)$, $i=1,2$. 
\end{Proof}
\subsubsection*{Ordered drifts} 
Theorem~\ref{pp4} also admits the following version
in the case of ordered drifts. 
\begin{corollary}
	\label{pp5}
        {\em (Increasing convex order)}. 
	 Assume that $X_0^{(1)} \leq X_0^{(2)}$ and
	$$
	\mu_1(t,x)\leq \mu_2(t,x)
	\quad
	\mbox{and}
	\quad 
	0 < \sigma_1(t,x)\leq \sigma_2(t,x),
	\qquad
	t\in [0,T], 
	\quad
	x \in \real,
	$$
	together with the following conditions: 
 \begin{itemize} 
\item[{\em (\hypertarget{C1}$C_1$)}] 
$\displaystyle g_1(t,x,y,z) \leq g_2(t,x,y,z)$, $t\in [0,T]$, $x,y \in \real$, $z\in \real_+$,
\item[{\em (\hypertarget{C2}$C_2$)}] 
$\displaystyle g_i(t,x,y,z)$ is non-decreasing in $z$ for $i=1$ or $i=2$, $x,y \in \real$, $z \in \real_+$, $t\in [0,T]$, 
\item[{\em (\hypertarget{C3}$C_3$)}] 
$\displaystyle g_i(t,x,y,z)$ is non-decreasing in $x$ for $i =1,2$, $x,y \in \real$, $z \in \real_+$, $t \in [0,T]$, 
\item[{\em (\hypertarget{C4}$C_4$)}] 
 $\displaystyle	(x,y) \mapsto f_i(t,x,y,z)$ and $(y,z) \mapsto f_i(t,x,y,z)$ 
   are both convex respectively on $\real^2$
   and $\real\times \real_+$ for $i=1,2$, $x,y \in \real$,
   $z \in \real_+$, $t\in [0,T]$.  
		\end{itemize} 
	Then we have $X_T^{(1)}\le_{g_1,g_2}^{\rm iconv} X_T^{(2)}$, i.e.,
	$$
	\cE_{g_1}\big[ \phi\big(X_T^{(1)}\big)\big] \leq \cE_{g_2}\big[ \phi\big(X_T^{(2)}\big)\big],
	$$
	for all non-decreasing convex functions $\phi(x)$ satisfying \eqref{growth.phi}. 
\end{corollary}
\begin{Proof} 
  Under (\hyperlink{C3}{$C_3$}),
  since $\phi(x) $ and $ g_i(t,x,y,z),i=1,2$, are non-decreasing in $ x $,
	by Proposition~\ref{nondecreasing}
	the solutions $u_1(t,x) $ and $u_2(t,x)$ of \eqref{du1-1} are 
	non-decreasing in $x $ and, as in the proof of 
	Theorem~\ref{pp4}, one can take $z\ge0$ 
	since $\displaystyle
	\frac{\partial u_i}{\partial x} (t,x) \ge0$.
	Assuming that e.g. $g_1(t,x,y,z) $ is non-decreasing in $z$
        under (\hyperlink{C2}{$C_2$}), 
	then by $ z\sigma_1(t,x)  \leq  z\sigma_2(t,x)$,
	$(t,x) \in [0,T]\times \real$, $z \in \real_+$, and
        (\hyperlink{C1}{$C_1$}), we have 
	$$ g_1(t,x,y,z\sigma_1(t,x)) \leq g_1(t,x,y,z\sigma_2(t,x)) \leq g_2(t,x,y,z\sigma_2(t,x)). 
	$$
	Combining the above with the inequality
        $ z\mu_1(t,x) \leq  z\mu_2(t,x), (t,x) \in [0,T]\times \real$,
        $z \in \real_+$, one finds 
	$f_1(t,x,y,z) \leq f_2(t,x,y,z)$, and by
        Theorem~\ref{pp4} we conclude that 
	$\cE_{g_1}\big[ \phi(X_T^{(1))}\big]
	\leq  \cE_{g_2}\big[ \phi \big(X_T^{(2)} \big)\big]$  
	for all convex non-decreasing functions $\phi(x)$ satisfying \eqref{growth.phi}.
\end{Proof}
When the drift coefficients $\mu(t,x)=\mu_1(t,x)=\mu_2(t,x)$
are equal and $g_i(t,x,y,z)$ is independent of $z$, $i=1,2$, 
the following proposition can be proved
for the convex $g$-ordering similarly to Corollary~\ref{pp5}, by applying
Theorem~\ref{pp3} which deals with convex ordering,
instead of Theorem~\ref{pp4}.
\begin{corollary} \label{pp6}
  {\em (Convex order with equal drifts)}.  
  Assume that $X_0^{(1)} = X_0^{(2)}$ and
	$$
	\mu_1(t,x)=\mu_2(t, x),\quad \mbox{and} \quad 0 < \sigma_1(t,x)\leq \sigma_2(t,x),
	\qquad
	t\in [0,T],
	\quad
	x \in \real,
	$$
 together with the conditions 
\begin{itemize} 
\item[{\em (\hypertarget{C'1}$C'_1$)}] 
   $	\displaystyle  g_i(t,x,y,z) = g_i(t,x,y)$
   is independent of $z\in \real$ for $i=1,2$, $t\in [0,T]$, $x,y \in \real$,
\item[{\em (\hypertarget{C'2}$C'_2$)}] 
   $g_1(t,x,y) \leq g_2(t,x,y)$, $t\in [0,T]$, $x,y \in \real$,
\item[{\em (\hypertarget{C'3}$C'_3$)}] 
$  	\displaystyle
   (x,y) \mapsto f_i(t,x,y,z)$ and $(y,z) \mapsto f_i(t,x,y,z)$
    are convex on $\real^2$ for $i=1,2$.  
		\end{itemize} 
	Then we have $X_T^{(1)}\le_{g_1,g_2}^{\rm conv} X_T^{(2)}$,  i.e., 
	$$
	\cE_{g_1} \big[\phi\big(X_T^{(1)}\big)\big] \leq \cE_{g_2} \big[ \phi\big(X_T^{(2)}\big)\big],
	$$
	for all convex functions $\phi(x)$ satisfying \eqref{growth.phi}. 
\end{corollary} 
We note that the convexity of $ u_1(t,x) $ and $ u_2(t,x) $
is not needed in the proof of Theorem~\ref{pp3}
when $\sigma_1(t,x) = \sigma_2(t,x)$,
and in this case we can remove Condition~(\hyperlink{B'2}{$B'_2$}) 
in Theorem~\ref{pp4} as in the next corollary.
\begin{corollary}
	\label{pp7} 
	{\em (Monotonic order with equal volatilities)}.  
        Assume that $X_0^{(1)} \leq X_0^{(2)}$ and $$ 0 < \sigma(t,x): = \sigma_1(t,x)=\sigma_2(t,x), \quad t\in [0,T], \quad x\in \real,
	$$ 
	together with the following conditions: 
 \begin{itemize} 
\item[{\em (\hypertarget{D1}$D_1$)}] 
$\displaystyle  \mu_1(t, x)\leq \mu_2(t, x)$, $x  \in \mathbb{R}$, $t\in [0,T]$,
\item[{\em (\hypertarget{D2}$D_2$)}] 
$g_1(t, x,y,z) \leq g_2(t, x,y,z)$  for all $(x,y,z)
	\in \real^2\times \real_+$, $t\in [0,T]$, 
\item[{\em (\hypertarget{D3}$D_3$)}] 
  $g_i(t, x,y,z)$ is non-decreasing in $x$ for $i=1,2$ and
  $(t,y,z) \in [0,T]\times \real \times \real_+$.
 \end{itemize} 
	Then we have
	$X_T^{(1)} \le_{g_1,g_2}^{\rm mon} X_T^{(2)}$, i.e., 
		$$
	\cE_{g_1} \big[\phi\big(X_T^{(1)}\big)\big] \leq \cE_{g_2} \big[ \phi\big(X_T^{(2)}\big)\big],
	$$
	for all non-decreasing functions $\phi(x)$ satisfying \eqref{growth.phi}. 
\end{corollary}
\begin{Proof} 
  Similarly to the proof of Corollary~\ref{pp4.1}, 
  under the condition $ \sigma_1(t,x) =\sigma_2(t,x)$ 
	the convexity of $ u_i(t,x)$ and the non-decreasing property
	of $ g_i(t,x,y,z)$ with respect to $ z $, $ i=1 $ or $ i=2 $,
	are no longer required.
	In addition, the condition 
	$$
        f_1(t,x,y,z) \leq f_2(t,x,y,z), \quad x,y \in \real, \ z\in \real_+, \ t\in [0,T], 
	$$
	clearly holds 
        from  (\hyperlink{D1}{$D_1$})-(\hyperlink{D2}{$D_2$}),
        and we can conclude as in the proof of Corollary~\ref{pp5}. 
\end{Proof}
             \section{Ordering with partially convex drifts} 
                        \label{s3.1}
                               Theorems~\ref{pp3} and \ref{pp4}
        require the convexity assumptions 
        (\hyperlink{B2}{$B_2$}) and (\hyperlink{B'2}{$B'_2$}) on 
        $$
        (x,y,z) \mapsto
        f_i(t,x,y,z): = z \mu_i(t,x)+ g_i(t,x,y,z\sigma_i(t,x))
        $$
        in $(x,y)$ and $(y,z)$ to hold for both $i=1,2$.
        In this section, we develop different 
        convex $g$-ordering results under weaker convexity
        conditions, based on a measurable function $\zeta (t,x)$ such that 
                         $$
\E \left[ \exp\left(\frac{1}{2}\int_0^T \left(
  \frac{\mu_i\big(t,X_t^{(i)}\big) - \zeta \big( t,X_t^{(i)} \big)}{\sigma_i \big(t,X_t^{(i)}\big)}
  \right)^2dt\right) \right] < \infty,
\quad i =1,2.
$$
As in Section~\ref{s3}, 
the proofs of Theorems~\ref{pp1}-\ref{pp2}
are first stated for ${\cal C}_b^3$ BSDE coefficients
as in Theorem~\ref{classical},   
        and then extended under (\hyperlink{A1}{$A_1$})-(\hyperlink{A4}{$A_4$}) 
using Proposition~\ref{dependence}. 
\begin{theorem}
\label{pp1}
{\em (Convex order)}.  
Assume that $X_0^{(1)}= X_0^{(2)}$ and
\begin{equation}
  \label{*101} 
  0 < \sigma_1(t,x) \leq \sigma_2(t,x), \qquad
  t\in [0,T], 
  \quad
  x \in \real,
\end{equation} 
  together with the conditions  
\begin{itemize} 
\item[{\em (\hypertarget{E1}$E_1$)}] 
$  \displaystyle  
\ f_1(t,x,y,z)
    \leq
  z \zeta (t,x) 
  \leq
 f_2(t,x,y,z)$, 
$ t\in [0,T]$, $ x,y,z \in \real$,
\item[{\em (\hypertarget{E2}$E_2$)}] 
$\   	\displaystyle
  (x,y) \mapsto f_i(t,x,y,z)$  and $(y,z) \mapsto f_i(t,x,y,z)$
  are convex on $\real^2$ for $i=1$ or $i=2$.  
\end{itemize} 
 Then we have $X_T^{(1)}\le_{g_1,g_2}^{\rm conv} X_T^{(2)}$,  i.e.,
\begin{equation}
\label{ceg} 
 \cE_{g_1}\big[\phi\big(X_T^{(1)}\big)\big] \leq \cE_{g_2}\big[\phi\big(X_T^{(2)}\big)\big],
\end{equation}
for all convex functions  $\phi(x)$ satisfying \eqref{growth.phi}. 
\end{theorem} 
\begin{Proof} 
  \noindent
  $(i)$ We start by assuming that the function $\phi$ and the coefficients
$\mu_i(t,\cdot)$,
$\sigma_i(t,\cdot)$
and $g_i(t,\cdot,\cdot,\cdot)$ are 
${\cal C}_b^3$ functions for all $t\in [0,T]$
as in Theorem~\ref{classical}, and that 
(\hyperlink{E2}{$E_2$}) holds with $i=1$.
Let
$$
 \theta_2 ( t , x) : = \frac{\mu_2(t,x) - \zeta ( t,x)}{\sigma_2 (t,x)},
 \qquad
 x\in \real,
 \quad
 t\in [0,T].
$$ 
 By the Girsanov theorem, the process 
$$
 \widetilde{B}_t := B_t +\int_0^{t}\theta_2 \big(s,X_s^{(2)}\big)  ds,
 \qquad
 t\in [0,T], 
$$ 
 is a standard Brownian motion under the probability measure
 $\mathbb{Q}_2$ defined by
$$
 \frac{d\mathbb{Q}_2}{d\mathbb{P}} := 
   \exp\left(-\int_0^T \theta_2 \big(s,X_s^{(2)}\big) dB_s -\frac{1}{2}\int_0^T \big( \theta_2 \big(s,X_s^{(2)}\big)\big)^2ds\right), 
$$ 
   and the forward SDEs \eqref{d3}-\eqref{d4} can be rewritten as 
$$
 \left\{
 \begin{array}{lll}
   dX^{(1)}_t & =
 \big(
 \mu_1 \big(t,X^{(1)}_t\big)-\theta_2  \big(t,X^{(2)}_t\big)\sigma_1 \big(t,X^{(1)}_t\big)\big)
 dt+\sigma_1 \big(t,X^{(1)}_t\big) d\widetilde{B}_t,
 & X_0^{(1)} = x_0^{(1)}, 
 \\
 \\
 dX^{(2)}_t & = 
 \zeta \big(t,X^{(2)}_t\big) dt
 + \sigma_2 \big(t,X^{(2)}_t\big) d\widetilde{B}_t,  & X_0^{(2)} = x_0^{(2)},
 \end{array}
 \right.
 $$ 
with the associated BSDEs 
$$
\left\{
\begin{array}{ll}
\nonumber
dY_t^{(1)} = - \big( g_1\big(t,X_t^{(1)},Y_t^{(1)},Z_t^{(1)}\big)+
Z_t^{(1)}
\theta_2 \big( t,X_t^{(2)}\big)
\big) dt + Z_t^{(1)}d\widetilde{B}_t, & Y_T^{(1)} =\phi\big( X_T^{(1)}\big),
  \\
  \\
  \nonumber
  dY_t^{(2)} = -\big( g_2\big(t,X_t^{(2)},Y_t^{(2)},Z_t^{(2)}\big)+
  Z^{(2)}_t \theta_2 \big( t,X_t^{(2)}\big)
  \big) dt + Z_t^{(2)} d\widetilde{B}_t, & Y_T^{(2)} =\phi\big( X_T^{(2)}\big).  
\end{array}
\right.
$$
      By Theorem~\ref{classical}
      we have
$Y^{(1)}_t = u_1 \big(t,X^{(1)}_t\big)$
and
$Y^{(2)}_t = u_2 \big(t,X^{(2)}_t\big)$,
where the functions
$u_1(t,x)$ and $u_2(t,x)$ are in ${\cal C}^{1,2}([0,T]\times\real)$
and solve the PDEs
		\begin{equation}
                  \label{du1} 
                  \frac{\partial u_i}{\partial t}(t,x)
                                    + \frac{1}{2}\sigma^2_i(t,x) \frac{\partial^2 u_i}{\partial x^2}(t,x)
                  + f_i\Big(t,x,u_i(t,x), \frac{\partial u_i}{\partial x}(t,x) \Big)=0,
		\end{equation}
                with $ u_i(T,x) =\phi(x) $, $i=1,2$.
                Applying  It\^o's formula to $ u_1\big(t, X_t^{(2)}\big)$ 
 and using \eqref{du1},    we have
                	\begin{eqnarray} 
                \nonumber
                u_1\big(t,X^{(2)}_t\big)
                & = & 
                u_1\big(0,X^{(2)}_0\big) 
                + \int_0^t 
                \frac{\partial u_1}{\partial s} \big(s,X^{(2)}_s\big)ds 
                + \int_0^t 
                                \zeta \big(s,X^{(2)}_s\big)
                               \frac{\partial u_1}{\partial x} \big(s,X^{(2)}_s\big)ds 
                \\
                \nonumber
               & & 
                + \frac{1}{2}\int_0^t \sigma_2^2\big(s,X^{(2)}_s\big) \frac{\partial^2 u_1}{\partial x^2}\big(s,X^{(2)}_s\big)ds 
                 + \int_0^t \sigma_2\big(s,X^{(2)}_s\big)\frac{\partial u_1}{\partial x} \big(s,X^{(2)}_s\big)d\widetilde{B}_s
                \\
                \nonumber
                & = & 
                u_1\big(0,X^{(2)}_0\big) 
                + \int_0^t 
                                \zeta \big(s,X^{(2)}_s\big)
                               \frac{\partial u_1}{\partial x} \big(s,X^{(2)}_s\big)ds 
                \\
                \nonumber 
              &  &  - \int_0^t f_1\Big(s,X^{(2)}_s,u_1\big(s,X^{(2)}_s\big),
                \frac{\partial u_1}{\partial x}\big(s,X^{(2)}_s\big) \Big)ds 
                \\
                \nonumber
               & & 
                + \frac{1}{2}\int_0^t \left(\sigma_2^2\big(s,X^{(2)}_s\big)- \sigma^2_1\big(s,X^{(2)}_s\big)\right)\frac{\partial^2 u_1}{\partial x^2}\big(s,X^{(2)}_s\big)ds 
                \\
                \nonumber
                & & + \int_0^t \sigma_2\big(s,X^{(2)}_s\big)\frac{\partial u_1}{\partial x} \big(s,X^{(2)}_s\big)d\widetilde{B}_s. 
                \end{eqnarray} 
                        Taking expectation at time $t=T $ under $\mathbb{Q}_2$,
                        we find 
                			\begin{eqnarray} 
                                          \nonumber 
			                    \E_{\mathbb{Q}_2} \big[ \phi \big(X^{(2)}_T\big) \big] 
                                          & = &
                                          u_1\big(0,X^{(2)}_0\big) 
                                          + \E_{\mathbb{Q}_2} \left[ \int_0^T
                                \zeta \big(s,X^{(2)}_s\big)
                                                 \frac{\partial u_1}{\partial x} \big(s,X^{(2)}_s\big)ds\right]
                          \nonumber \\
			& &  
			- \E_{\mathbb{Q}_2} \left[\int_0^T f_1\Big(s,X^{(2)}_s,u_1\big(s,X^{(2)}_s\big), \frac{\partial u_1}{\partial x}\big(s,X^{(2)}_s\big) \Big)ds\right] \nonumber \\
			& & + \frac{1}{2}\E_{\mathbb{Q}_2} \left[\int_0^T \left(\sigma_2^2\big(s,X^{(2)}_s\big)- 
			  \sigma^2_1\big(s,X^{(2)}_s\big)\right)\frac{\partial^2 u_1}{\partial x^2} \big(s,X^{(2)}_s\big)ds\right]. 
                        \nonumber
			\end{eqnarray} 
                                        Next, applying similarly It\^{o}'s formula to $ u_2\big(t, X_t^{(2)}\big)$ 
 and then taking expectation at $ t=T $ under
 $\mathbb{Q}_2$ we obtain, from \eqref{du1}, 	
			\begin{eqnarray} 
			  \E_{\mathbb{Q}_2} \big[ \phi \big(X^{(2)}_T\big) \big]
			  & = & u_2\big(0,X^{(2)}_0\big) 
                          - \E_{\mathbb{Q}_2} \left[\int_0^T f_2\Big(
			s,X^{(2)}_s,u_2\big(s,X^{(2)}_s\big), \frac{\partial u_2}{\partial x}\big(s,X^{(2)}_s\big) \Big)ds\right] \nonumber \\
			  & & + \E_{\mathbb{Q}_2} \left[\int_0^T \zeta \big(s,X_s^{(2)}\big) \frac{\partial u_2}{\partial x}\big(s,X^{(2)}_s\big)ds\right].
                         \nonumber 
			\end{eqnarray} 
 From Assumption~(\hyperlink{E1}{$E_1$}) and Condition~\eqref{*101} 
 we get 
			\begin{eqnarray} 
                          \label{weget}
                          \lefteqn{
      \! \! \! \! \!
                             u_2 \big(0,X^{(2)}_0\big)
-
u_1 \big(0,X^{(2)}_0\big)
 =
                          \E_{\mathbb{Q}_2} \left[\int_0^T 
			                          f_2\Big(s,X^{(2)}_s,u_2\big(s,X^{(2)}_s\big), \frac{\partial u_2}{\partial x}\big(s,X^{(2)}_s\big) \Big)ds \right] 
                          }
                            \\
                          \nonumber 
                          &	&
                          - \E_{\mathbb{Q}_2} \left[\int_0^Tf_1\Big(s,X^{(2)}_s,u_1\big(s,X^{(2)}_s\big), \frac{\partial u_1}{\partial x}\big(s,X^{(2)}_s\big) \Big)ds
			    \right]
                          \\
                          \nonumber 
			  & & - \E_{\mathbb{Q}_2} \left[\int_0^T \zeta \big(s,X_s^{(2)}\big) \frac{\partial u_2}{\partial x}\big(s,X^{(2)}_s\big)ds\right]
                          + \E_{\mathbb{Q}_2} \left[\int_0^T
                                                       \zeta \big(s,X_s^{(2)}\big)
                            \frac{\partial u_1}{\partial x} \big(s,X^{(2)}_s\big)ds\right] 
                          \\
\nonumber 
			& & + \frac{1}{2}\E_{\mathbb{Q}_2} \left[\int_0^T 
			\left(\sigma_2^2\big(s,X^{(2)}_s\big)- \sigma^2_1\big(s,X^{(2)}_s\big)
			\right)\frac{\partial^2 u_1}{\partial x^2} \big(s,X^{(2)}_s\big)ds\right] \\
                          \nonumber
			  & \ge &
                           \E_{\mathbb{Q}_2} \left[\int_0^T 
                             \left(
			     f_2\Big(s,X^{(2)}_s,u_2\big(s,X^{(2)}_s\big), \frac{\partial u_2}{\partial x}\big(s,X^{(2)}_s\big) \Big)
                             - \zeta \big(s,X_s^{(2)}\big)\frac{\partial u_2}{\partial x}\big(s,X^{(2)}_s\big) 
                             \right) ds \right]
                           \\
                             \nonumber 
			     &	& - \E_{\mathbb{Q}_2} \left[\int_0^T
                               \left(
                               f_1\Big(s,X^{(2)}_s,u_1\big(s,X^{(2)}_s\big), \frac{\partial u_1}{\partial x}\big(s,X^{(2)}_s\big) \Big)
                               - \zeta \big(s,X_s^{(2)}\big)
                               \frac{\partial u_1}{\partial x}\big(s,X^{(2)}_s\big)
                               \right) ds
			     \right]
                           \\
                           \nonumber 
                            & \geq & 0, 
			\end{eqnarray} 
                        where we have used
                        (\hyperlink{E1}{$E_1$}) and the fact that 
$\displaystyle \frac{\partial^2 u_1}{\partial x^2}(t,x) \ge 0$,
as follows from Theorem~\ref{convex}.  
 
\noindent
$(ii)$ The case $i=2$ in Assumption~(\hyperlink{E2}{$E_2$}) is dealt with
                        similarly by applying It\^{o}'s formula to 
                        $ u_2\big(t, X_t^{(1)}\big)$ and then to $ u_1\big(t, X_t^{(1)}\big)$,
                     and by taking expectation at $ t=T $ under the probability measure $\mathbb{Q}_1$ defined by
$$
 \frac{d\mathbb{Q}_1}{d\mathbb{P}} := 
   \exp\left(-\int_0^T \theta_1 \big(s,X_s^{(1)}\big) dB_s -\frac{1}{2}\int_0^T \big( \theta_1 \big(s,X_s^{(1)}\big)\big)^2ds\right), 
$$ 
   where
   $$
 \theta_1 ( t , x) : = \frac{\mu_1(t,x) - \zeta ( t,x)}{\sigma_1 (t,x)},
 \qquad
 x\in \real,
 \quad
 t\in [0,T]. 
 $$ 
 In this case, from
 \eqref{*101} and
 (\hyperlink{E1}{$E_1$}) 
 we get 
			\begin{eqnarray} 
                          \nonumber
                          \lefteqn{
	u_2 \big(0,X^{(2)}_0\big)
-
u_1 \big(0,X^{(2)}_0\big)
= 
\E_{\mathbb{Q}_1} \left[\int_0^T 
			                          f_2\Big(s,X^{(1)}_s,u_2\big(s,X^{(1)}_s\big), \frac{\partial u_2}{\partial x}\big(s,X^{(1)}_s\big) \Big)ds \right] 
                          }
                          \\
                          \nonumber 
                          			&	& - \E_{\mathbb{Q}_1} \left[\int_0^Tf_1\Big(s,X^{(1)}_s,u_1\big(s,X^{(1)}_s\big), \frac{\partial u_1}{\partial x}\big(s,X^{(1)}_s\big) \Big)ds
			    \right]
                          \\
                          \nonumber 
			  & &
                          - \E_{\mathbb{Q}_1} \left[\int_0^T
                            \zeta\big(s,X^{(1)}_s\big)
                             \frac{\partial u_2}{\partial x} \big(s,X^{(1)}_s\big)ds\right] 
+ \E_{\mathbb{Q}_1} \left[\int_0^T \zeta \big(s,X_s^{(1)}\big) \frac{\partial u_1}{\partial x}\big(s,X^{(1)}_s\big)ds\right]
                          \\
\nonumber 
			& & + \frac{1}{2}\E_{\mathbb{Q}_1} \left[\int_0^T 
			\left(\sigma_2^2\big(s,X^{(1)}_s\big)- \sigma^2_1\big(s,X^{(1)}_s\big)
			\right)\frac{\partial^2 u_2}{\partial x^2} \big(s,X^{(1)}_s\big)ds\right]                      \\
                           \nonumber 
                            & \geq & 0, 
			\end{eqnarray} 
since 
$\displaystyle \frac{\partial^2 u_2}{\partial x^2}(t,x) \ge 0$  
by Theorem~\ref{convex}.
By the relations
$Y_0^{(1)} = u_1 \big(0,X^{(1)}_0\big)$, 
$Y_0^{(2)} = u_2 \big(0,X^{(2)}_0\big)$ 
and
$X^{(1)}_0=X^{(2)}_0$
we conclude to $Y^{(2)}_0-Y^{(1)}_0 \geq 0$,
which shows \eqref{ceg}. 
The extension of \eqref{ceg} to coefficients satisfying
(\hyperlink{A1}{$A_1$})-(\hyperlink{A4}{$A_4$})
follows as in the proof of Theorem~\ref{pp3}.
\end{Proof}
The next proposition deals with the increasing convex order, for which only
the Conditions (\hyperlink{E'1}{$E'_1$})-(\hyperlink{E'2}{$E'_2$})
 and $X_0^{(1)} \leq X_0^{(2)}$ 
are required in addition to Condition~\eqref{sgm} and 
(\hyperlink{E'3}{$E'_3$}) below. 
\begin{theorem}
\label{pp2}
      {\em (Increasing convex order)}.
      Assume that
  $X_0^{(1)} \leq X_0^{(2)}$
  and
  \begin{equation}
    \label{sgm} 
  0 < \sigma_1(t,x) \leq \sigma_2(t,x), \quad
  t\in [0,T], \quad
  x \in \real,
\end{equation} 
  together with the conditions  
\begin{itemize} 
\item[{\em (\hypertarget{E'1}$E'_1$)}] 
$  \displaystyle  
\ f_1(t,x,y,z)
    \leq
  z \zeta (t,x) 
  \leq
 f_2(t,x,y,z)$, 
 $t\in [0,T]$, $x,y \in \real$, 
 $z\in \real_+$,
\item[{\em (\hypertarget{E'2}$E'_2$)}] 
$   	\displaystyle
  (x,y) \mapsto f_i(t,x,y,z)$ and $(y,z) \mapsto f_i(t,x,y,z)$
  are respectively convex on $\real^2$ and $\real\times \real_+$
  for $i=1$ or $i=2$, 
\item[{\em (\hypertarget{E'3}$E'_3$)}] 
$   	\displaystyle
   x \mapsto g_i(t,x,y,z)$ is non-decreasing on $\real$ for $i=1,2$, $ y\in \real, \, z \in \real_+, \, t \in [0,T] $.
\end{itemize}  
 Then we have $X_T^{(1)}\le_{g_1,g_2}^{\rm iconv} X_T^{(2)}$,  i.e.,
 $$
 \cE_{g_1}\big[\phi\big(X_T^{(1)}\big)\big] \leq \cE_{g_2}\big[\phi\big(X_T^{(2)}\big)\big],
 $$ 
        for all non-decreasing convex functions $\phi(x)$ satisfying \eqref{growth.phi}. 
\end{theorem}
\begin{Proof} 
  As in the proof of Theorem~\ref{pp1} 
  we start with ${\cal C}_b^3$ coefficients, and then 
  extend the conclusion to coefficients satisfying
 (\hyperlink{A1}{$A_1$})-(\hyperlink{A4}{$A_4$}) using Proposition~\ref{dependence}. 
  If $ \phi(x) $ and $ g_i(t,x,y,z)$ are non-decreasing in $x$ 
  by (\hyperlink{E'3}{$E'_3$}), $i=1,2$, then
  by Proposition~\ref{nondecreasing} the solutions $u_1(t,x) $ and $u_2(t,x) $
  of the PDE~\eqref{du1} are nondecreasing in $x$ and satisfy 
  $$
  \frac{\partial u_1}{\partial x} (t,X^{(2)}_t) \ge 0 \
  \quad
  \mbox{and}
  \quad
  \frac{\partial u_2}{\partial x} (t,X^{(2)}_t) \ge 0,
  $$
  a.s.,
  $t\in [0,T]$, 
  hence Conditions (\hyperlink{E1}{$E_1$})-(\hyperlink{E2}{$E_2$}) 
  only need to hold for $z\geq 0$,
  showing the sufficiency of
  (\hyperlink{E'1}{$E'_i$}),
  $i=1,2,3$. In addition,
  we have $u_1\big(0,X_0^{(1)}\big)=Y_0^{(1)} \leq u_1\big(0,X_0^{(2)}\big)$
  by the assumption $X_0^{(1)} \leq X_0^{(2)}$, hence 
  by repeating arguments in the proof of Theorem~\ref{pp1}
  for $i=1$ we find 
 by \eqref{weget} that 
    \begin{eqnarray} 
          \nonumber
          \lefteqn{
  	Y^{(2)}_0-Y^{(1)}_0
  	\ge u_2\big(0,X_0^{(2)}\big) - u_1\big(0,X_0^{(2)}\big)
          }
  \\
  \nonumber 
    & =	&          
        \E_{Q_2} \left[\int_0^T 
  	f_2\Big(s,X^{(2)}_s,u_2\big(s,X^{(2)}_s\big), \frac{\partial u_2}{\partial x}\big(s,X^{(2)}_s\big) \Big)ds \right] 
  \\
  \nonumber 
    &	& - \E_{Q_2} \left[\int_0^Tf_1\Big(s,X^{(2)}_s,u_1\big(s,X^{(2)}_s\big), \frac{\partial u_1}{\partial x}\big(s,X^{(2)}_s\big) \Big)ds
    \right]
  \\
  \nonumber 
  & & - \E_{Q_2} \left[\int_0^T \zeta \big(s,X_s^{(2)}\big) \frac{\partial u_2}{\partial x}\big(s,X^{(2)}_s\big)ds\right]
  + \E_{Q_2} \left[\int_0^T
  \zeta\big(s,X_s^{(2)}\big)
  \frac{\partial u_1}{\partial x} \big(s,X^{(2)}_s\big)ds\right] 
  \\
    \nonumber 
& & + \frac{1}{2}\E_{Q_2} \left[\int_0^T 
  \left(\sigma_2^2\big(s,X^{(2)}_s\big)- \sigma^2_1\big(s,X^{(2)}_s\big)
  \right)\frac{\partial^2 u_1}{\partial x^2} \big(s,X^{(2)}_s\big)ds\right]
  \\
  \nonumber 
  & \ge & 0, 
 \end{eqnarray}
        under Assumption~(\hyperlink{E'2}{$E'_2$}) for $i=1$. 
        The case $i=2$ is treated similarly according to the proof of
        Theorem~\ref{pp1}. 
\end{Proof}
\section{Comparison in $g$-risk measures} \label{s5}
A $g$-risk measure is a mapping $\rho : L^2(\Omega, \cF_T) \to \R $
satisfying the following conditions.
\begin{definition} Let $g$ satisfy Conditions~(\hyperlink{A2}{$A_2$})-(\hyperlink{A3}{$A_3$}) and $\xi \in  L^2(\Omega, \cF_T)$. 
\begin{enumerate}[1)]
\item The static $g$-risk measure
        is defined in terms of $g$-evaluation as
	\begin{equation}
	  \nonumber
          \rho^g(X):= \cE_g [ -\xi ]. 
	\end{equation}
      \item The dynamic $g$-risk measure
        is defined in terms of conditional $g$-evaluation as
		\begin{equation}
		  \nonumber
                  \rho^g_t(\xi):= \cE_g [ -\xi \mid \cF_t ],
                  \qquad t\in[0,T]. 
		\end{equation}
\end{enumerate}
\end{definition}
We refer to \cite{gianin2006risk} for 
the relations between coherent and convex risk measures,
and the $g$-expectation.

\medskip

\noindent
 We note that, taking $ \phi(x):=-x$
in \eqref{co2}, $ X^{(1)} \le_{g_1,g_2}^{conv} X^{(2)} $
                implies $ \rho^{g_1}\big(X^{(1)}\big) \leq \rho^{g_2}\big(X^{(2)}\big)$.
                In addition, we have
                $\rho^{g_1}\big( \phi\big(X_T^{(1)} \big)\big) \leq \rho^{g_2}\big( \phi\big(X_T^{(2)} \big) \big)$ for all convex function $ \phi(x) $ if and only if
$$X_T^{(2)} \le_{g^{(-1)}_1,g^{(-1)}_2}^{\rm conv} X_T^{(1)}, 
$$
				where
                                $$
                                g^{(-1)}_1(t,x,y,z):=-g_1(t,x,-y,-z)
                                \quad and \quad
                                g^{(-1)}_2(t,x,y,z)=-g_2(t,x,-y,-z), 
                                $$ 
 as from Lemma~\ref{lem} with
                 $ a=-1$ we have 
                  $\cE_{g^{(-1)}}[\phi(X_T)]
                	= -\cE_{g^{(1)}} [ -\phi(X_T)]$. 
A stochastic ordering via $G$-expectations
        has also been defined in \cite{tian2016uncertainty}
 by combining \eqref{co2} with the inequality 
		\begin{equation}
                \nonumber 
                  -\rho^{g_1} \big(\phi\big(X^{(1)}\big)\big)
                  =
                  -\cE_{g_1}\big[ -\phi\big(X^{(1)}\big)\big] \le
                  -\cE_{g_2}\big[ - \phi\big( X^{(2)}\big)\big]
                  =-\rho^{g_2} \big( \phi\big(X^{(2)}\big)\big), 
                \end{equation}
 where 
$-\cE_{g_i}\big[-\phi\big(X_T^{(i)}\big)\big]$
and
$\cE_{g_i}\big[\phi\big(X_T^{(i)}\big)\big]$
respectively represent bid and ask prices of the contingent claim in financial markets, $i=1,2$.

\medskip

\noindent 
Theorems~\ref{pp3} and \ref{pp1} admit
the following versions for the comparison
of risks. First, we have the next consequence of
Theorem~\ref{pp3} and Lemma~\ref{lem} below,
where we let
$$
    f_i^{(-1)} ( t,x,y,z): = 
      	z\mu_i (t,x)
  	+
  	g^{(-1)}_i (t,x,y,z \sigma_1(t,x)),
\quad 
t\in \real_+, \ x,y,z\in \real, \ i = 1,2.
$$
 
\begin{corollary}
	          \label{pp9}
	Assume that
	$X_0^{(1)}= X_0^{(2)}$
	and
	$$0 < \sigma_1(t,x)\leq \sigma_2(t,x),
	\qquad 	t\in [0,T], \quad x\in \real,
	$$
	together with the conditions 
         \begin{itemize} 
         \item[{\em (\hypertarget{F1}$F_1$)}] 
$
	\displaystyle  
  	f^{(-1)}_1(t,x,y,z)
  	\leq
  	f^{(-1)}_2(t,x,y,z)
$, $t\in [0,T]$, $x,y,z \in \real$,
         \item[{\em (\hypertarget{F2}$F_2$)}] 
  	$ \displaystyle
        (x,y) \mapsto f_i^{(-1)}(t,x,y,z)$
        and $(y,z) \mapsto f_i^{(-1)}(t,x,y,z)$
        are convex on $\real^2$ for $i=1,2$.
\end{itemize}  
	Then we have
	$$
	-\cE_{g_1}\big[-\phi\big(X_T^{(1)}\big)\big] \leq -\cE_{g_2}\big[-\phi\big(X_T^{(2)}\big)\big], 
	$$ 
	for all convex functions $\phi(x)$ satisfying \eqref{growth.phi}. 
\end{corollary}
Similarly, we have the next consequence
 of Theorem~\ref{pp1} and Lemma~\ref{lem}. 
\begin{corollary} 
    \label{pp10}
  Assume that
  $X_0^{(1)}= X_0^{(2)}$
  and
  $$0 < \sigma_1(t,x)\leq \sigma_2(t,x),
  \qquad 	t\in [0,T],\quad x\in \real,
  $$
 together with the conditions 
\begin{itemize} 
   \item[{\em (\hypertarget{G1}$G_1$)}] 
$	\displaystyle  
  	f^{(-1)}_1(t,x,y,z)
  \leq
  z \zeta (t,x) 
    	\leq
  	f^{(-1)}_2(t,x,y,z)$, 
  $t\in [0,T]$, $x,y,z \in \real$,
\item[{\em (\hypertarget{G2}$G_2$)}] 
$   	\displaystyle
   (x,y) \mapsto f_i^{(-1)}(t,x,y,z)$
   and $(y,z) \mapsto f_i^{(-1)}(t,x,y,z)$
   are convex on $\real^2$ for $i=1$ or $i=2$. 
  	\end{itemize} 
		Then we have
                $$
                -\cE_{g_1}\big[-\phi\big(X_T^{(1)}\big)\big] \leq -\cE_{g_2}\big[-\phi\big(X_T^{(2)}\big)\big], 
                $$ 
                for all convex functions $\phi(x)$ satisfying \eqref{growth.phi}. 
\end{corollary} 
We note that Corollaries~\ref{pp9} and \ref{pp10} can be applied 
to Example~\ref{ex6.1} 
below for the comparison of $g_i$-risk measures
with $g_i(t,x,y,z)$ linear in $y$ and $z$,
in which case $g_i^{(-1)}(t,x,y,z)=g_i(t,x,y,z)$ and 
the bid and ask prices 
$-\cE_{g_i}\big[-\phi\big(X_T^{(i)}\big)\big] = \cE_{g_i}\big[\phi\big(X_T^{(i)}\big)\big]$
are equal, $i=1,2$.

\medskip
 
Furthermore, we can also derive results for the increasing convex and
monotonic orderings of $g$-risk measures under Conditions~(\hyperlink{F1}{$F_1$})-(\hyperlink{F2}{$F_2$})
and
(\hyperlink{G1}{$G_1$})-(\hyperlink{G2}{$G_2$}). 
For example, if
(\hyperlink{F1}{$F_1$})-(\hyperlink{F2}{$F_2$})
or 
(\hyperlink{G1}{$G_1$})-(\hyperlink{G2}{$G_2$})  
 only holds for $ z \ge 0$ and $ g_i$ is non-decreasing in $ x $, $i=1,2 $, we then get
versions of Corollaries~\ref{pp9}-\ref{pp10}
for the increasing convex $g$-risk comparisons
as in Theorems~\ref{pp4} and \ref{pp2}. 
Similarly,
under additional the assumption
$ \sigma_1(t,x) = \sigma_2(t,x)$
and by removing
(\hyperlink{F2}{$F_2$}) and
(\hyperlink{G2}{$G_2$}), we can obtain versions of
Corollaries~\ref{pp9}-\ref{pp10}
for the monotonic $g$-risk ordering as in
Corollary~\ref{pp4.1}. 

\section{Application examples}
\label{s4}
In the following examples we
\textcolor{newcolor}{ aim
 at comparing option prices
 of the form $Y_0^{(i)}:=\cE_{g_i}\big[\phi\big(X_T^{(i)}\big)\big]$}
for two risky assets
      with positive prices 
      $\big(X_t^{(i)}\big)_{t\in [0,T]}$, $i=1,2$, given by
      \begin{equation}
      \nonumber
      dX_t^{(i)}=X_t^{(i)} a_i\big(t,X_t^{(i)}\big)dt + X_t^{(i)} b_i\big( t, X_t^{(i)} \big) dB_t,
      \qquad
      i=1,2, 
      \end{equation}
      where the coefficients
      $\mu_i (t,x) = x a_i(t,x)$ and
      $\sigma_i (t,x) = x b_i(t,x)$
      satisfy
      (\hyperlink{A1}{$A_1$}),
      $i=1,2$, for example $a_i(t,x)$ and $b_i(t,x)$
      can be bounded functions.
\textcolor{newcolor}{      
  Example~\ref{ex6.1} compares option prices
  in classical expectation
  for standard self-financing portfolios
      with $ X_t^{(1)}=X_t^{(2)}:=X_t$, $a_1(t,x) =a_2(t,x) :=r$,
      $\sigma_t = b_1(t,X_t)$,
      with a misspecified volatility coefficient
      $b_2\big(t,x\big)$ such that $\sigma_t \leq b_2(t,X_t)$, $a.s.$, 
      and is consistent with Theorem~6.2 in \cite{elkjs}.
      }
\begin{example}
        \label{ex6.1}
 Taking 
\begin{equation}
  \label{gi}
  g_i(t,x,y,z) := -ry - z
      \frac{a_i( t,x )-r}{b_i( t,x)}, \qquad i=1,2,
  \end{equation} 
under the conditions
\begin{equation}
  \label{under}
  X_0^{(1)} = X_0^{(2)}
      \qquad
      \mbox{and}
      \qquad 0< b_1(t,x) \leq b_2(t,x),
\quad t\in [0,T], \ x>0, 
\end{equation} 
      we have
      $X_T^{(1)} \le_{g_1,g_2}^{\rm conv} X_T^{(2)} $,
      i.e. the values of the self-financing portfolios
      hedging the claim payoffs 
      $\phi\big(X_T^{(1)}\big)$ and $\phi\big(X_T^{(2)}\big)$
      satisfy 
      $$ 
      \cE_{g_1} \big[ \phi\big(X_T^{(1)}\big)\big] \leq \cE_{g_2}\big[ \phi\big(X_T^{(2)}\big)\big], 
      $$ 
      for all convex functions $\phi(x)$ satisfying \eqref{growth.phi}.
      \end{example}
      \begin{Proof} 
 Consider the risk-free asset $E_t := E_0 e^{rt}$ and
        the portfolio valued
      $$
      V^{(i)}_t : = p^{(i)}_t X^{(i)}_t + q^{(i)}_t E_t,
      \qquad t\in \real_+,
      $$
      where $p^{(i)}_t$ is the quantity
      of risky assets and $q^{(i)}_t$ is the quantity of risk-free assets. 
      When the strategy $(p^{(i)}_t, q^{(i)}_t)_{t\in \real_+}$ is self-financing, 
      we have 
      \begin{align}
      dV_t^{(i)} &  =q^{(i)}_tdE_t + p^{(i)}_tdX_t^{(i)}
      \\
      & =\big(
      rV_t^{(i)} +
      \theta_i\big( t,X_t^{(i)}\big)
      p^{(i)}_t X_t^{(i)} b_i \big( t,X_t^{(i)}\big) \big)
      dt +
      p^{(i)}_tX_t^{(i)} b_i \big( t,X_t^{(i)}\big)dB_t, 
      \nonumber 
      \end{align}
 where 
      $$
      \theta_i( t,x)
      := \frac{a_i( t,x )-r}{b_i( t,x)},
      \qquad
      i=1,2. 
      $$
 Hence, letting
      $$
      Z_t^{(i)}:=p^{(i)}_tX_t^{(i)} b_i \big( t,X_t^{(i)}\big),
      $$
      and discounting as 
      $$
      \widetilde{V}_t^{(i)} :=e^{-rt} V_t^{(i)},
      \quad 
      \widetilde{X}_t^{(i)} :=e^{-rt} X_t^{(i)},
      \quad 
      \mbox{and}
      \quad 
      \widetilde{Z}_t^{(i)} :=e^{-rt} Z_t^{(i)},
      $$
      with $V_T^{(i)}=\phi\big(X_T^{(i)}\big)$,
      we find the linear BSDE 
      \begin{equation}
      \nonumber
      \widetilde{V}_t^{(i)}=\widetilde{V}_T^{(i)}
      + \int_t^T \tilde{g}_i\big(s,\widetilde{X}_s^{(i)}, \widetilde{V}_s^{(i)}, \widetilde{Z}_s^{(i)}\big) ds
      - \int_t^T \widetilde{Z}_s^{(i)}dB_s,
      \end{equation}
      with $\tilde{g}_i(t,x,y,z) = - z \theta_i(t,xe^{rt})$,
      $i=1,2$. 
 Since 
      $$
      xz ( a_i (t,xe^{rt}) - r ) + \tilde{g}_i(t,x,y,z x b_i(t,xe^{rt}) ) 
      = 0,
      \quad 
      x,y,z \in (0,\infty ) \times \real^2, \ t\in [0,T],
            $$
      $i=1,2$, we check that Conditions~(\hyperlink{B1}{$B_1$}) and (\hyperlink{E1}{$E_1$})
      are satisfied (with $\zeta (t,x) = 0$) 
      together with
      (\hyperlink{B2}{$B_2$}) and (\hyperlink{E2}{$E_2$}), 
      hence, under \eqref{under},  
      Theorems~\ref{pp3} and \ref{pp1} show that 
            $$
      \widetilde{V}_0^{(1)} = \cE_{\widetilde{g}_1} \big[ e^{-rT}\phi \big(e^{rT} \widetilde{X}_T^{(1)}\big)\big] \leq  \widetilde{V}_0^{(2)} = \cE_{\widetilde{g}_2}\big[ e^{-rT}\phi\big(e^{rT}\widetilde{X}_T^{(2)}\big)\big],
      $$ 
      that is $\widetilde{X}_T^{(1)} \le_{\widetilde{g}_1,\tilde{g}_2}^{\rm conv} \widetilde{X}_T^{(2)} $. 
      \textcolor{newcolor}{Therefore the price of the first claim is 
        upper bounded by the price of the second claim, i.e. } 
      $$ 
      \cE_{g_1} \big[ \phi\big(X_T^{(1)}\big)\big] \leq \cE_{g_2}\big[ \phi\big(X_T^{(2)}\big)\big], 
      $$ 
      for all convex functions $\phi(x)$ satisfying \eqref{growth.phi},
      or $X_T^{(1)} \le_{g_1,g_2}^{\rm conv} X_T^{(2)} $, 
      with $g_i(t,x,y,z)$ as in \eqref{gi}, 
      $i=1,2$.
      \end{Proof}
      \textcolor{newcolor}{
        The next example considers the kernel
        $g (t,x,y,z) := \alpha(t)|z|$ with $\alpha(t) >0$, $t\in [0,T]$,
        and shows that the second risky asset 
        would be preferred over the first one 
        by risk-seeking investors whose preferences are modeled by $g$ 
        since in this case $\cE_g [ \ \! \cdot \ \! ]$
        is represented as in \eqref{g_alpha-a},
        see also Theorem~1.37 in \cite{sriboonchita2009stochastic}.
         \begin{example}
Consider the kernel 
      	\begin{equation} 
\nonumber 
          g_i(t,x,y,z)=
     	  g (t,x,y,z) :=\alpha(t)|z|, \quad z \in \real,
          \quad i=1,2,
      	\end{equation} 
        where $\alpha(t)$ is a positive bounded function on $[0,T]$.
 Then, assuming that 
      \begin{align}
      X_0^{(1)} \leq X_0^{(2)},\quad a_1(t,x) \leq a_2(t,x), \quad \mbox{and} \quad b_1(t,x) \leq b_2(t,x),   \ t\in [0,T], \ x>0, \label{ex6.1b_cond}
      \end{align}
 and
 that the functions $x\mapsto x a_i(t, x)$ and $x\mapsto x b_i(t, x)$ 
 are convex in $x\in \real_+$ for $i=1,2$ and $t\in [0,T]$, 
      we have
      \begin{equation}
        \label{fjkdslfd} 
      \cE_g \big[ \phi \big(X_T^{(1)}\big)\big] \leq \cE_g \big[ \phi \big(X_T^{(2)}\big)\big],
      \end{equation}
      for all non-decreasing convex functions $\phi(x)$ satisfying \eqref{growth.phi}. 
              \end{example}
         \begin{Proof}
         Under \eqref{ex6.1b_cond}, we check that
         when $\alpha(t)>0$, $t\in [0,T]$, we have 
         $$ xza_1(t,x) + z \alpha(t)b_1(t,x) \leq xza_2(t,x) + z \alpha(t)b_2(t,x),
         \quad \ x, z \in \real_+,
         $$
         and we conclude by Theorem~\ref{pp4}.
               \end{Proof}
         In the case $\alpha(t) <0$, $t\in [0,T]$,
         the function $ z\mapsto 
         g (t,x,y,z) =\alpha(t)|z|$ is concave on $\real$, 
         and a similar result can be obtained
         for risk-averse investors
         by taking $b_1(t,x) \ge b_2(t,x)$
         in \eqref{ex6.1b_cond}, provided that 
         the functions 
         and $x\mapsto x a_i(t, x)$, $x\mapsto x b_i(t, x)$ are concave in $x\in \real_+$ 
         for $=1,2$ and $t\in [0,T]$, 
         and that the function $\phi$ in \eqref{fjkdslfd} is non-decreasing concave. 
      }
       
      \medskip
       
      \noindent
      \textcolor{newcolor}{
         In Examples~\ref{ex6.2}-\ref{ex6.4}
      we consider portfolios under constraints,
      in which case the BSDE generators are sublinear functions.}
      In Example~\ref{ex6.2} we consider the comparison of option
prices for a standard self-financing portfolio
and a second self-financing hedging 
portfolio in which borrowing occurs at the rate $R \geq r $, 
as in Example~1.1 in \cite{el1997backward}.
\begin{example}
 \label{ex6.2}
 Taking
      $$
      g_1(t,x,y,z) := -ry -z\frac{a_1( t,x )-r}{b_1( t,x)}
      $$
      and
            $$
      g_2(t,x,y,z) := -ry
      -z\frac{a_2( t,x )-r}{b_2( t,x)} +(R-r)\left( y- \frac{z}{b_2(t,x)} \right)^-,
      $$ 
      where $w^- = - \min (w,0)$, and under the conditions
      \begin{equation}
      \label{a12} 
      X_0^{(1)} = X_0^{(2)}
      \qquad
      \mbox{and}
      \qquad
      0< b_1(t,x) \leq b_2(t,x),
      \quad 
      t \in [0,T], \ x > 0, 
      \end{equation} we have
      $X_T^{(1)} \le_{g_1,g_2}^{\rm conv} X_T^{(2)}$,
      i.e. 
      $$ 
      \cE_{g_1} \big[ \phi \big(X_T^{(1)}\big) \big] \leq \cE_{g_2} \big[ \phi \big(X_T^{(2)}\big) \big], 
      $$ 
      for all convex functions $\phi(x)$ satisfying \eqref{growth.phi}, i.e., the fair price of the unconstrained portfolio is less than
      that of the one with constraints. 
      \end{example}
\begin{Proof} 
      The first portfolio value is the discounted wealth process
of Example~\ref{ex6.1}, which satisfies the BSDE 
      \begin{equation}
     \nonumber
     \widetilde{V}^{(1)}_t= e^{-rT}\phi\big(e^{rT}\widetilde{X}_T^{(1)}\big)
     + \int_t^T
     \tilde{g}_1(s,\widetilde{X}_s^{(1)}, \widetilde{V}^{(1)}_s, \widetilde{Z}_s^{(1)}\big) ds - \int_t^T \widetilde{Z}_s^{(1)}dB_s,
     \end{equation}
      where $\theta_1(t,x) := ( a_1(t,x)-r ) / b_1(t,x)$,
      with the generator $\tilde{g}_1(t,x,y,z) :=-z\theta_1(t,xe^{rt})$.
      In the second portfolio the investor is only
      allowed to borrow at the rate $ R \geq r $,
      which yields the discounted wealth process
      \begin{equation}
      \nonumber
      \widetilde{V}^{(2)}_t= e^{-rT}\phi\big(e^{rT}\widetilde{X}_T^{(2)}\big)
      + \int_t^T \tilde{g}_2\big( s,\widetilde{X}_s^{(2)},\widetilde{V}^{(2)}_s, \widetilde{Z}_s^{(2)} \big) ds
            - \int_t^T \widetilde{Z}_s^{(2)}dB_s,
      \end{equation}
      which is a BSDE with the generator
      $$
      \tilde{g}_2(t,x,y,z)
      : = -z\theta_2(t,xe^{rt}) +(R-r)\left( y- \frac{z}{b_2(t,xe^{rt})} \right)^-.
      $$ 
      We check that
      \begin{align} 
        xz ( a_1(t,xe^{rt})  - r ) + \tilde{g}_1(t,x,y,z x b_1(t,xe^{rt}) )
                 & = 0
        \\
        &\leq (R-r)(y-xz)^-
        \nonumber \\ 
        & = z x ( a_2(t,xe^{rt}) - r ) + \tilde{g}_2(t,x,y,z xb_2(t,xe^{rt}) ),
        \nonumber
      \end{align} 
      $x,y,z \in (0,\infty ) \times \real^2$, $t\in [0,T]$,
 where both functions $(x,y) \mapsto (R-r)(y-xz)^-$
      and $(y,z) \mapsto (R-r)(y-xz)^-$
      are convex,
      hence
      (\hyperlink{B2}{$B_2$}) and
      (\hyperlink{E2}{$E_2$}) are satisfied.
       In addition,
       (\hyperlink{B1}{$B_1$}) and (\hyperlink{E1}{$E_1$}) hold (with $\zeta(t,x)=0$)  
      hence by \eqref{a12}, Theorems~\ref{pp3} and \ref{pp1} both show that
 $$
      \widetilde{V}^{(1)}_0 = \cE_{\widetilde{g}_1} \big[ e^{-rT}\phi\big(e^{rT}\widetilde{X}_T^{(1)}\big)\big] \le
      \widetilde{V}^{(2)}_0 = \cE_{\widetilde{g}_2}\big[ e^{-rT}\phi\big(e^{rT}\widetilde{X}_T^{(2)}\big)\big],
      $$ 
      or
      $$\cE_{g_1} \big[ \phi\big(X_T^{(1)}\big)\big] \leq \cE_{g_2}\big[ \phi\big(X_T^{(2)}\big)\big],
      $$
      with $g_i(t,x,y,z) = -ry + \tilde{g}_i(t,xe^{-rt},y,z)$, $i=1,2$,
      for all convex functions $\phi(x)$ satisfying \eqref{growth.phi}, which shows the part (a).
\end{Proof} 
      \textcolor{newcolor}{
        In particular, when  $ X_t^{(1)}=X_t^{(2)}:=X_t$,
        $a_1(t,x) =a_2(t,x) :=a(t,x)$
        and
        $b_1(t,x) = b_2\big(t,x\big):=b(t,x)$, Example~\ref{ex6.2} shows that
        $$ 
       \cE_{g_1} \big[ \phi \big(X_T\big) \big] \leq \cE_{g_2} \big[ \phi \big(X_T\big) \big], 
       $$ 
       for all convex functions $\phi(x)$ satisfying \eqref{growth.phi}.
       Here,
       using the same underlying risky asset $(X_t)_{t \in [0, T]}$,
       the first agent hedges the contingent claim $\phi(X_T)$ by
       a self-financing strategy without
       borrowing money,
       while the second agent
       hedges the same claim
       under constraints on the
       difference between the borrowing and lending rates.
       In this case,
       the initial investment of the
       second agent is higher as we have $V_0^{(2)} \ge V_0^{(1)}$. 
      }
      
      \medskip 

      \noindent
       In the next example we
assume that both self-financing hedging portfolios
require borrowing at the rate $R \geq r$. 
     \begin{example}
 Taking 
      $$g_i(t,x,y,z) := -ry
      -z\frac{a_i( t,x )-r}{b_i( t,x)} +(R-r)\left( y- \frac{z}{b_i (t,x)} \right)^-, \qquad i=1,2,
      $$ 
      under the conditions
         \begin{equation}
        \label{and} 
        X_0^{(1)} = X_0^{(2)}
        \qquad
        \mbox{and} 
        \qquad 0< b_1(t,x) \leq b_2(t,x),
\quad 
            t \in [0,T], \ x > 0, 
      \end{equation} 
      we have
      $X_T^{(1)} \le_{g_1,g_2}^{\rm conv} X_T^{(2)} $,
      i.e. 
      $$ 
      \cE_{g_1} \big[ \phi\big(X_T^{(1)}\big)\big] \leq \cE_{g_2}\big[ \phi\big(X_T^{(2)}\big)\big], 
      $$ 
      for all convex functions $\phi(x)$ satisfying \eqref{growth.phi}. 
      \end{example}
      \begin{Proof} 
        We consider 
        two portfolios constructed as in Example~\ref{ex6.2}, 
     with discounted wealth processes given by the BSDEs
      \begin{align*}
      \nonumber
      \widetilde{V}_t^{(i)} = &\,\,  e^{-rT}\phi\big(e^{rT}\widetilde{X}_T^{(i)}\big)
      + \int_t^T       \tilde{g}_i\big( s,\widetilde{X}_s^{(i)},\widetilde{V}_s^{(i)}, \widetilde{Z}_s^{(i)}\big)ds
         - \int_t^T \widetilde{Z}_s^{(i)}dB_s, \quad i=1,2, 
      \end{align*}
      with the generators
      $$
      \tilde{g}_i(t,x,y,z)
      : = -z\theta_i(t,xe^{rt}) +(R-r)\left( y-\frac{z}{b_i(t,xe^{rt})} \right)^-, \quad i=1,2. 
      $$ 
    We check that 
      \begin{align} 
      xz ( a_1(t,xe^{rt})  - r ) + \tilde{g}_1(t,x,y,z x b_1(t,xe^{rt}) ) 
      & = (R-r)(y-xz)^-
      \nonumber \\ 
      & = z x ( a_2(t,e^{rt}x) - r ) + \tilde{g}_2(t,x,y,z xb_2(t,xe^{rt}) ),
      \nonumber
      \end{align} 
      $x,y,z \in (0,\infty ) \times \real^2$, $t\in [0,T]$,
      hence
      (\hyperlink{B1}{$B_1$}) is satisfied.
      Here, Theorem~\ref{pp1} cannot be applied because
      (\hyperlink{E1}{$E_1$}) is
           not satisfied.
           However, (\hyperlink{B2}{$B_2$}) is satisfied since 
      both functions $(x,y) \mapsto (R-r)(y-xz)^-$
      and $(y,z) \mapsto (R-r)(y-xz)^-$
      are convex, hence by \eqref{and}, Theorem~\ref{pp3} shows that
      $$
      \widetilde{V}_0^{(1)} =
      \cE_{\widetilde{g}_1} \big[ e^{-rT}\phi\big(e^{rT}\widetilde{X}_T^{(1)}\big)\big] \le
      \widetilde{V}_0^{(2)} = \cE_{\widetilde{g}_2}\big[ e^{-rT}\phi\big(e^{rT}\widetilde{X}_T^{(2)}\big)\big],
      $$ 
      or
      $$
      \cE_{g_1} \big[ \phi\big(X_T^{(1)}\big)\big] \le
      \cE_{g_2}\big[ \phi\big(X_T^{(2)}\big)\big],
      $$
      with $g_i(t,x,y,z) = -ry + \tilde{g}_i(t,xe^{-rt},y,z)$, $i=1,2$,
      for all convex functions $\phi(x)$ satisfying \eqref{growth.phi}.
      \end{Proof}
\noindent       
      The next Example~\ref{ex6.4} is based on three risky asset prices
      $\big(X_t^{(1)}\big)_{t\in [0,T]}$,
      $\big(X_t^{(2)}\big)_{t\in [0,T]}$ and $\big(X_t^{(3)}\big)_{t\in [0,T]}$,
      see \S3.2 of \cite{jouini1995arbitrage}.
      The portfolio of the first investor is based on the
      risky asset $X_t^{(1)}$ and on
      the risk-free asset $E_t := E_0e^{rt}$
      as in Example~\ref{ex6.1}. 
      On the other hand, the second investor
            is longing $X_t^{(2)}$ and $ E_t $ while 
      shorting $X_t^{(3)} $ and
the       risk-free asset 
$\overline{E}_t := \overline{E}_0e^{ R t}$,
with $R \geq r$.
 \begin{example}
              \label{ex6.4}
              In addition to $\big(X_t^{(1)}\big)_{t\in [0,T]}$
              and $\big(X_t^{(2)}\big)_{t\in [0,T]}$,
              consider a third asset with positive price
                            $\big(X_t^{(3)}\big)_{t\in [0,T]}$
      given by 
      \begin{equation}
      \nonumber
      dX_t^{(3)}=X_t^{(3)} a_3\big(t,X_t^{(2)}\big)dt
      + X_t^{(3)} b_3\big( t, X_t^{(2)} \big)
      dB_t, \quad t\in[0,T]. 
      \end{equation}
Let 
                    $$
                    g_1(t,x,y,z) := -ry - z \frac{a_1(t,x) - r}{b_1 (t,x)},
                                        $$
 and                    
 $$
 g_2(t,x,y,z) :=
 -ry
 -z^+\frac{a_2( t,x )-r}{b_2( t,x)} +z^-\frac{a_3( t,x )-r}{b_3( t,x)}
 +(R -r)\left(
      y-\frac{z^+}{b_2(t,x)}+\frac{z^-}{b_3(t,x)}\right)^-, 
      $$
 and assume that 
      \begin{equation}
        \label{thta} 
        \theta_2(t,x):=\frac{a_2(t,x)-r}{b_2(t,x)}
 \leq 
        \theta_3(t,x) := \frac{a_3(t,x)- r}{b_3(t,x)}, 
      \qquad t\in [0,T], \ x>0. 
\end{equation} 
Then, under the conditions  
\begin{equation}
  \label{under2} 
                        X_0^{(1)} = X_0^{(2)}
                      \qquad
                      \mbox{and}
                      \qquad 
                      0< b_1(t,x) \leq b_2(t,x), \quad t \in [0,T], \, x > 0,
\end{equation} 
 we have 
$X_T^{(1)} \le_{g_1,g_2}^{\rm conv} X_T^{(2)} $, i.e.       
      $$ 
      \cE_{g_1} \big[ \phi\big(X_T^{(1)}\big)\big] \leq \cE_{g_2}\big[ \phi\big(X_T^{(2)}\big)\big],
      $$ 
 for all convex functions $\phi(x)$ satisfying \eqref{growth.phi}, i.e., the fair price of the short-selling constrained portfolio is greater than
 that of the one without constraints.  
\end{example} 
      \begin{Proof}
                            Under Condition~\eqref{thta} 
                    the model is without arbitrage
                    by Theorems~3.1 and 3.2
                    in \cite{jouini1995arbitrage}.
                    In the optimal solution of
                    Theorem~3.2 therein, 
      the optimal hedging strategy at time $t$ for the second investor 
      is to long
      $$  \Delta^{(2)}_{t}
      :=\left(\frac{\partial v}{\partial x}\big(t,X_t^{(2)}\big)\right)^+
      $$
      units of $ X_t^{(2)} $, and to short 
      $$
      \Delta^{(3)}_{t}:=-\frac{b_2\big(t,X_t^{(2)}\big)X_t^{(2)}}{b_3\big(t,X_t^{(2)}\big)X_t^{(3)}}\left(\frac{\partial v}{\partial x}\big(t,X_t^{(2)}\big)\right)^-,
      $$
      units of $X_t^{(3)}$,
      while longing $\big(\Delta_t^{(0)}\big)^+$ units of $ E_t$,
     and shorting $-\big(\Delta_t^{(0)}\big)^-  $ units of $\overline{E}_t$,
      where
$$ 
      \Delta_t^{(0)}:=v\big(t,X_t^{(2)}\big)-X_t^{(2)} \Delta^{(2)}_t - X_t^{(3)} \Delta^{(3)}_t,
      \qquad
      t\in [0,T].
 $$ 
      In other words, the discounted portfolio
asset price processes
 $      \widetilde{V}_t^{(i)} :=e^{-rt} V_t^{(i)}$ and  
$      \widetilde{X}_t^{(i)} :=e^{-rt} X_t^{(i)}$ 
       satisfy the BSDEs 
      \begin{equation}
      \nonumber
      \widetilde{V}_t^{(i)}= e^{-rT}\phi\big(e^{rT}\widetilde{X}_T^{(i)}\big)
      + \int_t^T
      \tilde{g}_i\big( s,\widetilde{X}_s^{(i)}, \widetilde{V}_s^{(i)},\widetilde{Z}_s^{(i)}\big) ds
      - \int_t^T \widetilde{Z}_s^{(i)}dB_s, \quad
      i=1,2,
      \end{equation}
     with 
        $$
     \widetilde{Z}_t^{(2)}:=
     e^{-rt} X_t^{(2)} b_2\big(t,X_t^{(2)}\big)\frac{\partial v}{\partial x}\big(t,X_t^{(2)}\big),
      \qquad t \in [0,T], 
      $$
and the generators
      $\tilde{g}_1(t,x,y,z) := - z \theta_1(t,xe^{rt})$
and
$$     \tilde{g}_2(t,x,y,z) := -z\theta_2(t,xe^{rt}) + z^-\big(\theta_3(t,xe^{rt})-\theta_2(t,xe^{rt} )\big) +(R -r)\left(
      y-\frac{z^+}{b_2(t,xe^{rt})}+\frac{z^-}{b_3(t,xe^{rt})}\right)^-.
$$ 
 Hence       the second portfolio price $
      V_t^{(2)} = v\big(t,X_t^{(2)}\big)$ satisfies
                    the PDE
      \begin{align*}
        \frac{\partial v}{\partial t}(t,x) &+rx\frac{\partial v}{\partial x}(t,x) +\frac{1}{2}x^2b_2(t,x)\frac{\partial^2 v}{\partial x^2}(t,x)
        + x \left( \theta_3(t,x)-\theta_2(t,x)\right)
        b_2(t,x)\left( \frac{\partial v}{\partial x}(t,x)\right)^-  \\
      & -rv(t,x)+(R -r)\left( v(t,x)-x\left(\frac{\partial v}{\partial x}(t,x)\right)^+ +x\frac{b_2(t,x)}{b_3(t,x)}\left(\frac{\partial v}{\partial x}(t,x)\right)^- \right)^- =0. 
      \end{align*}
      Then, by Conditions~\eqref{thta}-\eqref{under2} we have
      \begin{eqnarray*}
        \lefteqn{
          xz ( a_1(t,xe^{rt}) - r ) + \tilde{g}_1(t,x,y,z x b_1(t,xe^{rt}) ) 
          = 0 }
        \\
        &\leq &
        xz^-b_2(t,x) \big(\theta_3(t,xe^{rt})-\theta_2(t,xe^{rt} )\big)
 +(R -r)\left(y-xz^+ +\frac{b_2(t,xe^{rt})}{b_3(t,xe^{rt})}xz^-\right)^-  \\
     &  = & xz ( a_2(t,xe^{rt}) - r ) + \tilde{g}_2(t,x,y,z xb_2(t,xe^{rt}) ),
 \qquad
 (x,y,z) \in \real_+\times \real^2, \ t\in [0,T],
      \end{eqnarray*}
      hence (\hyperlink{E2}{$E_2$}) is satisfied.
      Condition~(\hyperlink{E1}{$E_1$})
      is satisfied with $\zeta (t,x) = 0$ 
      hence, under Condition~\eqref{under2},  
      Theorem~\ref{pp1} shows that 
            $$ 
      \widetilde{V}_0^{(1)} = \cE_{\widetilde{g}_1} \big[ e^{-rT}\phi \big(e^{rT} \widetilde{X}_T^{(1)}\big)\big] \leq \widetilde{V}_0^{(2)} = \cE_{\widetilde{g}_2}\big[ e^{-rT}\phi\big(e^{rT}\widetilde{X}_T^{(2)}\big)\big],
      $$
      or
      $$ 
      \cE_{g_1} \big[ \phi\big(X_T^{(1)}\big)\big] \leq \cE_{g_2}\big[ \phi\big(X_T^{(2)}\big)\big],
      $$ 
      with $g_i(t,x,y,z) = -ry + \tilde{g}_i(t,xe^{-rt},y,z)$, $i=1,2$,
            for all convex functions $\phi(x)$ satisfying \eqref{growth.phi},
      that is $X_T^{(1)} \le_{g_1,g_2}^{\rm conv} X_T^{(2)} $.           
      Note that here, Theorem~\ref{pp3} cannot be applied since the function 
      $$
      (x,y) \longmapsto xz ( a_2(t,xe^{rt}) - r ) + \tilde{g}_2(t,x,y,z xb_2(t,xe^{rt}) )
      $$
      may not be convex, hence (\hyperlink{B2}{$B_2$})
      is not satisfied. 
      \end{Proof}
\textcolor{newcolor}{
  We note that the conclusions of the above examples also imply
  the comparison of terminal portfolios values 
      $$
      \cE_g \big[ u\big(V_T^{(1)}\big)\big] \leq \cE_g \big[ u\big(V_T^{(2)}\big)\big],
      $$
      for all non-decreasing convex utility and payoff functions $u$ and $\phi$, 
      since $u\circ \phi(x)$ is also non-decreasing and convex,
      hence in this case
      the second portfolio would be preferred
      over the first portfolio by risk-seeking investors
        whose preferences are modeled by $g$. 
}
\section{Convexity of nonlinear PDE solutions} 
In this section, we extend the convexity result
Theorem~1.1 in \cite{bian2008convexity}
for nonlinear PDEs 
under a weaker convexity condition on the nonlinear drift 
$(x,y,z)\mapsto f(t,x,y,z)$ 
in the one-dimensional case, as required by applications
in finance, see the nonlinear Examples~\ref{ex6.2}-\ref{ex6.4}. 
For this we remark that, in our one-dimensional
setting, the constant rank
Theorem~2.3 in \cite{bian2008convexity},
see also Theorem~1.2 in \cite{bian2009},
only requires 
convexity of the nonlinear drift $f(t,x,y,z)$ 
    in $(x,y)\in \real^2$ for every
    $(t,z)\in [0,T]\times\real$,
    instead of global convexity in $(x,y,z)$.
    Precisely, we note that Condition~(2.6) in Theorem 2.3 of \cite{bian2008convexity}
    reduces to \eqref{Q.condition} below.
    \begin{theorem}(Theorem 2.3, \cite{bian2008convexity}). 
  \label{rank.theo}
  Assume that $u(t,x)$ is a ${\cal C}^{2,4}\big([0,T)\times \real\big)$
 convex solution of the PDE 
          \begin{equation}
          \label{gen.pde}
            \frac{\partial u}{\partial t}(t,x)+F\Big(t,x,u(t,x),\frac{\partial u}{\partial x}(t,x),\frac{\partial^2 u}{\partial x^2}(t,x)\Big)=0,
          \end{equation}
          and that
          $F(t,x,y,z,w)$ is a ${\cal C}^{1,2}\big([0,T)\times \real^4\big)$ function that
    satisfies the elliptic condition 
  \begin{equation}
    \nonumber 
    \frac{\partial F}{\partial w}\Big(t,x,u(t,x),\frac{\partial u}{\partial x}(t,x),\frac{\partial^2 u}{\partial x^2}(t,x)\Big) > 0,
    \quad t\in [0,T], \ x\in \real, 
\end{equation}
    and 
    \begin{align} 
      \nonumber
        \frac{\partial^2 F}{\partial x^2}\Big(t,x,u(t,x),\frac{\partial u}{\partial x}(t,x),0\Big) & + 2b\frac{\partial^2 F}{\partial x \partial y}\Big(t,x,u(t,x),\frac{\partial u}{\partial x}(t,x),0\Big)
      \\
      \label{Q.condition}
       &  + b^2\frac{\partial^2 F}{\partial y^2}\Big(t,x,u(t,x),\frac{\partial u}{\partial x}(t,x),0\Big)
      \ge 0, 
    \end{align} 
    for all
    $t\in [0,T]$, $x\in \real$, and $b \in \real$. 
Then the sign
$\displaystyle {\rm sgn} \left( \frac{\partial^2 u}{\partial x^2}(t,x) \right)$ of
$\displaystyle \frac{\partial^2 u}{\partial x^2}(t,x) $ is constant in $x\in \real$
        for any $ t\in(0,T) $, 
	and we have
        $$
        \displaystyle {\rm sgn} \left( \frac{\partial^2 u}{\partial x^2}(s,x) \right)
        \geq
            {\rm sgn} \left( \frac{\partial^2 u}{\partial x^2}(t,x) \right),
            \qquad
            0 \leq s \leq t < T.
            $$ 
\end{theorem}
\textcolor{newcolor}{
  We note that in our one-dimensional setting, the rank
  of $\displaystyle \frac{\partial^2 u}{\partial x^2}(t,x)$
  is $\{0,1\}$-valued, with the relation
  $\displaystyle {\rm sgn} \Big( \frac{\partial^2 u}{\partial x^2}(t,x) \Big)
  \equiv {\rm rank} \Big( \frac{\partial^2 u}{\partial x^2}(t,x) \Big)$ 
  provided that 
  $x\mapsto u(t,x)$ convex in $x \in\real$,
  $t \in [0,T]$. }
      By adapting arguments of \cite{bian2008convexity}, 
    using the constant rank Theorem~\ref{rank.theo}
      and a new Lemma~\ref{lemma}, 
    we will prove the following Theorem~\ref{convex} which
 has been used in the proofs of
Theorems~\ref{pp3}-\ref{pp4},
Corollaries~\ref{pp5}-\ref{pp6}, and 
Theorems~\ref{pp1}-\ref{pp2}
and Corollaries~\ref{pp9}-\ref{pp10}.

\medskip

\textcolor{newcolor}{
  We note that
  by Theorem~2.3 and Condition~(2.6) in \cite{bian2008convexity},
  convexity of the solution $u(t,x)$
  of the PDE \eqref{du12} is ensured
  by the joint convexity of
  $(x,y,z) \mapsto f(t,x,y,z):=F(t,x,y,z,0)$, 
  where $F(t,x,y,z,w)$ is the semilinear function
  $$F(t,x,y,z,w):=z\mu(t,x)
  +g\left(t,x,y,z\sigma(t,x)\right) +\frac{w}{2}\sigma^2(t,x).
  $$
  However, this joint convexity condition is too strong for our
  applications in mathematical finance,
  and in Theorem~\ref{convex}
  we show that it can be relaxed into the 
  partial convexity of $F(t,x,y,z,0)$ in $(x,y)$ and in $(y,z)$.
  Namely,} in Theorem~\ref{convex} we extend Theorem~1.1 of \cite{bian2008convexity} on the convexity of the solution $u(t,x)$ 
by only assuming convexity in $(x,y)$ and in $(y,z)$ of the function $f(t,x,y,z)$ in
\eqref{f.func} below, 
instead of joint convexity in $(x,y,z)\mapsto f(t,x,y,z)$.
  \textcolor{newcolor}{For this, we start by assuming
  additional regularity conditions on coefficients,
  and then use approximation arguments
  with help of the continuous dependence (stability)
  Theorem~\ref{dependence}, see Lemma~\ref{lemma}. 
}
\begin{theorem}
	\label{convex}
  Assume that the coefficients $ \mu $, $ \sigma$, 
  $ g$ and $\phi$ satisfy (\hyperlink{A1}{$A_1$})-(\hyperlink{A4}{$A_4$}). 
	Suppose that $u(t,x)$ is a ${\cal C}^{1,2}([0,T ) \times \mathbb{R})$
          solution of \eqref{du12} with terminal condition $ u(T,x)=\phi(x)$,  
		and that the function
\begin{equation}
\label{f.func} f(t,x,y,z) :=z\mu(t,x)  +g\left(t,x,y,z\sigma(t,x)\right)
\end{equation}
 satisfies the following conditions: 
\begin{itemize} 
\item[{\em (\hypertarget{H1}$H_1$)}] 
   $\displaystyle  (x,y) \mapsto f(t,x,y,z)$ is convex on $\real^2$ for every $(t,z)\in [0,T]\times \real$,
\item[{\em (\hypertarget{H2}$H_2$)}] 
      $(y,z) \mapsto f(t,x,y,z)$ is convex on $\real^2$ for every $(t,x)\in [0,T]\times\real$.
\end{itemize} 
	Then the function $x \mapsto u(t,x) $ is convex on $ \real$ for all $t \in [0,T]$,
	provided that $u(T,x)=\phi(x) $ convex in $x \in \real$. 
\end{theorem}
\begin{Proof}
	 We proceed by extending the proof argument of Theorem~1.1 in \cite{bian2008convexity}
 by using Theorem~\ref{rank.theo} 
 and an approximation argument.
 We start by assuming that the following conditions, 
 which are stronger than 
 	 (\hyperlink{A1}{$A_1$})-(\hyperlink{A4}{$A_4$}),
 hold for some $\eta \in (0,1)$.
    \begin{itemize} 
 \item[(\hypertarget{H3}{$H_3$})] 
   $\mu(\cdot,\cdot)$, $\sigma(\cdot,\cdot)$,
   $g(\cdot,\cdot,y,z) \in {\cal C}_b^{1+\eta / 2 ,2 + \eta }\big([0,T]\times \real \big) $
   for all $y,z\in \real$, and $ \phi(\cdot) \in C_b^{4+\eta } (\real) $, 
\item[(\hypertarget{H4}{$H_4$})] 
         $\sigma(\cdot,\cdot)$ satisfies the bound 
 	\begin{equation}
 	\label{cc2} 
 	0 < c   \leq \sigma(t,x) \leq C, \qquad t\in [0,T], \quad x \in \real,
 	\end{equation} 
 	for some constants $c,C>0$. 
\item[(\hypertarget{H5}{$H_5$})] 
 For some $ C>0$ and $\alpha > 0$ we have
 $$\displaystyle \Big|\frac{\partial^2 g}{\partial z^2}(t,x,y,z)\Big| \leq \frac{C}{(1+x^2)^{\alpha +1}},
 \qquad
 (t,x,y,z) \in [0,T]\times \real^3.
$$
\end{itemize} 
 Under Conditions~(\hyperlink{H3}{$H_3$})-(\hyperlink{H4}{$H_4$})  
  the function $ F(t,y,z,w) $ defined as 
\begin{equation}
  \nonumber 
  F(t,x,y,z,w):= f(t,x,y,z) + \frac{w}{2}\sigma^2(t,x), \quad t \in [0,T), \, x,\,y,\,z,\,w \in \real,
\end{equation}
is in ${\cal C}^{1,2}\big([0,T)\times \real^4\big)$
  and the solution $u(t,x)$ of \eqref{gen.pde} is in
  ${\cal C}^{2,4}([0,T)\times \real)$
    by Theorem~\ref{classical}.
    Besides, we note that by Theorem~\ref{classical2},
  the first and the second partial derivatives of $u(t,x)$ with respect to $x$
  are bounded uniformly in $(t,x)\in [0,T]\times \real$.
  Setting $ h(x) := (1+x^2)^{\alpha + 1}$, 
  for any $K\in \real$ and $\varepsilon >0$, we define  
  $$
  v_K(t,x):=e^{-Kt}h(x)
  \quad
  \mbox{and}
  \quad   u_{\varepsilon} (t,x) :=u(t,x)+\varepsilon v_K(t,x),
	\qquad x\in \real, \quad t\in [0,T]. 
	$$
	Next, we let 
	$$
	E_{\varepsilon} :=\left\{
	(t,x)\in [0,T]\times\real \ : \ \frac{\partial^2 u_{\varepsilon}}{\partial x^2} (t,x)\leq 0 \right\} 
	$$ 
	and suppose that $ E_{\varepsilon}  \ne \emptyset $. From the relation $h''(x) \ge (1+x^2)^{\alpha }$ and the
        bound $\displaystyle \left|\frac{\partial^2 u}{\partial x^2}(t,x)\right|\leq C$
        we get 
	$$\frac{\partial^2 u_{\varepsilon}}{\partial x^2}(t,x) \ge \varepsilon e^{-Kt}(1+x^2)^{\alpha }-C,
	$$
	therefore there exists $R_\varepsilon>0$ such that
        $\displaystyle \frac{\partial^2 u_{\varepsilon}}{\partial x^2}(t,x) > 0$ for all $|x| \ge R_\varepsilon$,
	and we have $ E_{\varepsilon}
	  \subseteq [0,T]\times B(0,R_\varepsilon)$,
	where $B(0,R_\varepsilon) $ is the centered open ball with radius $ R_\varepsilon $,
	so that $E_{\varepsilon} $ is compact. Consequently, the supremum 
	$$ \tau_0 :=\sup \{t \in [0,T] \ : \ (t,x)\in E_{\varepsilon}  \mbox{ for some } x \in \real \}
	$$
	is attained at some $(\tau_0,x_0) \in E_{\varepsilon} $ with $x_0 \in B(0,R_\varepsilon)$,
	such that $\displaystyle \frac{\partial^2 u_{\varepsilon}}{\partial x^2}(\tau_0,x_0) \leq 0 $.
	In addition, by the convexity assumption on $x\mapsto u(T,x)$ we have
	$$
	\frac{\partial^2 u_{\varepsilon}}{\partial x^2}(T,x)=\frac{\partial^2 u}{\partial x^2}(T,x)+\varepsilon \frac{\partial^2 v_K}{\partial x^2}(T,x) \ge \varepsilon e^{-KT} h''(x) >0,
        \qquad x\in \real, 
	$$
	hence $\tau_0 <T$ and by the continuity of $u_{\varepsilon}$ we have 
	$\displaystyle \frac{\partial^2 u_{\varepsilon}}{\partial x^2}(\tau_0,x) \ge 0 $,
	$x\in \real$, 
	since 
	$\displaystyle \frac{\partial^2 u_{\varepsilon}}{\partial x^2}(t,x) > 0 $ for all $t\in (\tau_0, T)$
	and $ x \in \real$.
        	Consequently, the function
	$ u_{\varepsilon}(t,x) $ is convex in $x$ on $[ \tau_0,T ]  \times B(0,R_\varepsilon)$.

\medskip
	
\noindent
On the other hand, we note that 
	$\displaystyle \frac{\partial^2 u_{\varepsilon}}{\partial x^2}(\tau_0,x_0) = 0$ 
		for $x_0 \in B(0,R_\varepsilon )$,
        and that $u_{\varepsilon}(t,x)$ satisfies the equation
	$$ \frac{\partial u_{\varepsilon}}{\partial t}(t,x)
	+ F_{K, \varepsilon } \left(t,x,u_{\varepsilon} (t,x),\frac{\partial u_{\varepsilon}}{\partial x}(t,x),\frac{\partial^2 u_{\varepsilon}}{\partial x^2}(t,x)\right)
	=0, $$
	where 
	\begin{eqnarray} 
          \label{fk} 
	F_{K, \varepsilon }(t,x,y,z,w) & := &  
	- \varepsilon \frac{ \partial v_K}{\partial t}(t, x) +\frac{1}{2}\sigma^2(t,x)
        \left( w-\varepsilon \frac{\partial^2 v_K}{\partial x^2}(t,x)\right)
        \\
        \nonumber 
       & & +f\left(t,x,y-\varepsilon v_K(t,x),z-\varepsilon \frac{\partial v_K}{\partial x}(t,x)
        \right).
\end{eqnarray} 
	By the constant rank Theorem~\ref{rank.theo} and Lemma~\ref{lemma} below, 
	we deduce that
        for small enough $T=T(\varepsilon , \alpha )>0$
        the second derivative  
        $\displaystyle \frac{\partial^2 u_{\varepsilon}}{\partial x^2} (t,x) $
	vanishes on $[ \tau_0 , T) \times B(0,R_\varepsilon )$ hence $\tau_0=T$,
	which is a contradiction showing that $ E_{\varepsilon} = \emptyset $.
	Therefore we have
	$$
	\frac{\partial^2 u}{\partial x^2}(t,x)+\varepsilon \frac{\partial^2 v_K}{\partial x^2}(t,x)
	=
	\frac{\partial^2 u_{\varepsilon}}{\partial x^2}(t,x) > 0,
        \qquad
        (t,x) \in [0,T]\times \real, 
	$$
	and after letting
	$\varepsilon$ tend to $0$, we conclude that
	$$
	\frac{\partial^2 u}{\partial x^2}(t,x) \geq 0,
        \qquad
        (t,x) \in [0,T]\times \real, 
	$$
        for small enough $T>0$. 
        This conclusion extends to all $T>0$ by decomposing $[0,T]$
        into subintervals of lengths at most $T(\varepsilon , \alpha )>0$.
        Finally, we relax the above
        Conditions~(\hyperlink{H3}{$H_3$})-(\hyperlink{H5}{$H_5$})
        under the hypotheses
        (\hyperlink{A1}{$A_1$})-(\hyperlink{A4}{$A_4$})
by applying the above argument to sequences  
$(\mu_n)_{n\geq 1}$,
$(\sigma_n)_{n\geq 1}$,
$(g_n)_{n\geq 1}$ of ${\cal C}_b^{2,3}$ functions and $(\phi)_{n\geq 1}$ of ${\cal C}_b^5$ 
 functions satisfying
 (\hyperlink{H3}{$H_3$}) 
 and (\hyperlink{A1}{$A_1$})-(\hyperlink{A4}{$A_4$}), 
 and converging pointwise respectively to
$\mu$, $\sigma$, $g$, and $ \phi $
while preserving the convexity of
the approximations $ (\phi_n)_{n\ge 1} $ and
$(f_n)_{n\geq 1}$
defined in \eqref{f.func}, as well as Condition~\eqref{cc2}. 
In order to satisfy (\hyperlink{H5}{$H_5$}), we replace $g_n$ with 
$\tilde{g}_n$ obtained by smoothing out the piecewise ${\cal C}^1$ function 
\begin{eqnarray*} 
x & \mapsto &  
\left(
g_n(t,-n,y,z)+ ( x+n) \frac{\partial g_n}{\partial x} (t,-n,y,z)
\right)
      {\bf 1}_{(-\infty, -n )} (x) 
\\
& &
+ g_n(t,x,y,z) {\bf 1}_{[-n,n]} (x) 
      + \left( g_n(t,n,y,z) + ( x-n)
      \frac{\partial g_n}{\partial x}(t,n,y,z) \right) {\bf 1}_{(n,\infty )} (x) 
      ,
\end{eqnarray*} 
     by convolution in $x$ with the Gaussian kernel $e^{-n x^2/2}/\sqrt{2\pi/n}$,  
and we conclude by the continuous dependence Proposition~\ref{dependence}.
\end{Proof}
The proof of Theorem~\ref{convex} relies on the following lemma.
\begin{lemma}
  \label{lemma}
  Under Conditions (\hyperlink{H1}{$H_1$})-(\hyperlink{H5}{$H_5$}) above, 
        for $T=T(\varepsilon , \alpha )>0$ small enough we can choose $K\in \real$ 
	  such that the function $(x,y) \mapsto F_{K , \varepsilon } (t,x,y,z,0) $ 
in \eqref{fk} satisfies Condition~\eqref{Q.condition}.
\end{lemma}
\begin{Proof}
  We need to show that
  \begin{eqnarray}
  \nonumber 
    S(b)  &:= & \frac{\partial^2 F_{K , \varepsilon }}{\partial x^2}\Big(t,x,u_{\varepsilon}(t,x),\frac{\partial u_{\varepsilon}}{\partial x}(t,x),0\Big) +2b\frac{\partial^2 F_{K , \varepsilon }}{\partial x \partial y}\Big(t,x,u_{\varepsilon}(t,x),\frac{\partial u_{\varepsilon}}{\partial x}(t,x),0\Big) \qquad 
    \\
    \nonumber
    & & +b^2\frac{\partial^2 F_{K , \varepsilon }}{\partial y^2}\Big(t,x,u_{\varepsilon}(t,x),\frac{\partial u_{\varepsilon}}{\partial x}(t,x),0\Big)
  \\
  \nonumber
   & \ge & 0,
  \end{eqnarray}
   for all  $ b \in \real $. 
 By \eqref{fk}, we have 
   \begin{align*}
   & S(b)  
   = 2\frac{\partial^2 f}{\partial x^2} -2\Big(\varepsilon\frac{\partial v_K}{\partial x}(t,x)-b \Big)\frac{\partial^2 f}{\partial x \partial y}+ \frac{1}{2}\Big(\varepsilon\frac{\partial v_K}{\partial x}(t,x)-b \Big)^2\frac{\partial^2 f}{\partial y^2} \\
   & +\frac{1}{2}\Big(\varepsilon\frac{\partial v_K}{\partial x}(t,x)-b \Big)^2\frac{\partial^2 f}{\partial y^2}+2\varepsilon\frac{\partial^2 v_K}{\partial x^2}(t,x)\Big(\varepsilon \frac{\partial v_K}{\partial x}(t,x) - b\Big) \frac{\partial^2 f}{\partial y \partial z} +2\Big(\varepsilon\frac{\partial^2 v_K}{\partial x^2}(t,x)\Big)^2 \frac{\partial^2 f}{\partial z^2 }  \\
   & -\frac{\partial^2 f}{\partial x^2}-\Big(\varepsilon\frac{\partial^2 v_K}{\partial x^2}(t,x) \Big)^2\frac{\partial^2 f}{\partial z^2}-\varepsilon\left( \frac{1 }{2}\frac{\partial^2 \sigma^2 }{\partial x^2}(t,x)+\frac{\partial f}{\partial y}+2\frac{\partial^2 f}{\partial x \partial z}\right) \frac{\partial^2 v_K}{\partial x^2}(t,x)\\
   &   
   -\varepsilon \frac{\partial^3 v_K}{\partial x^2\partial t}(t,x)-\varepsilon \left(\frac{\partial f}{\partial z}+2\sigma(t,x)\frac{\partial \sigma }{\partial x}(t,x)\right)\frac{\partial^3 v_K}{\partial x^3}(t,x) 
   - \frac{\varepsilon}{2}\sigma^2(t,x)\frac{\partial^4 v_K}{\partial x^4}(t,x), 
   \end{align*}
   for all $ b \in \real $, 
   where the derivatives of $f$ are evaluated at the point
   $\big(t,x,u_\varepsilon(t,x),\frac{\partial u_\varepsilon}{\partial x}(t,x) \big)$. 
Since $ f(t,x,y,z) $ is convex in $ (x,y) $, and in $ (y,z) $,  we have
	$$2\frac{\partial^2 f}{\partial x^2} -2\Big(\varepsilon\frac{\partial v_K}{\partial x}(t,x)-b \Big)\frac{\partial^2 f}{\partial x \partial y}+\frac{1}{2}\Big(\varepsilon\frac{\partial v_K}{\partial x}(t,x)-b \Big)^2\frac{\partial^2 f}{\partial y^2} \ge 0,
	$$
	and 
	$$ \frac{1}{2}\Big(\varepsilon\frac{\partial v_K}{\partial x}(t,x)-b \Big)^2\frac{\partial^2 f}{\partial y^2}+2\varepsilon\frac{\partial^2 v_K}{\partial x^2}(t,x)\Big(\varepsilon \frac{\partial v_K}{\partial x}(t,x) - b\Big) \frac{\partial^2 f}{\partial y \partial z} +2\Big(\varepsilon\frac{\partial^2 v_K}{\partial x^2}(t,x)\Big)^2 \frac{\partial^2 f}{\partial z^2 }   \ge 0,
        $$
	for all $b \in \real$.
        To conclude it suffices to 
        show the inequality 
	\begin{eqnarray*}
	S'& :=&-\frac{\partial^2 f}{\partial x^2}-\varepsilon\Big(\frac{\partial^2 v_K}{\partial x^2}(t,x) \Big)^2\frac{\partial^2 f}{\partial z^2} -\left( \frac{1}{2}\frac{\partial^2 \sigma^2 }{\partial x^2}(t,x) + \frac{\partial f}{\partial y}+2\frac{\partial^2 f}{\partial x \partial z}\right) \frac{\partial^2 v_K}{\partial x^2}(t,x)\\
	&   & 
	- \frac{\partial^3 v_K}{\partial x^2\partial t}(t,x)	- \left(\frac{\partial f}{\partial z}+2\sigma(t,x)\frac{\partial \sigma }{\partial x}(t,x)\right)\frac{\partial^3 v_K}{\partial x^3}(t,x) 	- \frac{1}{2}\sigma^2(t,x)\frac{\partial^4 v_K}{\partial x^4}(t,x) \\
         & \ge & 0,
	\end{eqnarray*}
        Thanks to Conditions~(\hyperlink{H3}{$H_3$})-(\hyperlink{H4}{$H_4$}), 
        we find that there exists $C >0$ such that
        $$
          \Big|\frac{\partial^2 f}{\partial x^2}(t,x,y,z)\Big| \leq C \left(
 1 + z^2 \Big|\frac{\partial^2 g}{\partial x^2}(t,x,y,z\sigma (t,x))\Big| 
\right), \qquad 
       \Big|\frac{\partial f}{\partial y}(t,x,y,z)\Big| \leq C,
       $$
       and
       $$
       \Big|\frac{\partial^2 f}{\partial x \partial z}(t,x,y,z)\Big|\leq
       C\left(
       1+|z|
       \Big|\frac{\partial^2 g}{\partial z^2}(t,x,y,z)\Big|
       \right), \qquad 
       \Big|\frac{\partial f}{\partial z}(t,x,y,z)\Big| \leq C(1+|x|), 
       $$
        for all $(t,x,y,z) \in [0,T]\times \real^3$, hence
        by Condition~(\hyperlink{H5}{$H_5$})
        for some $C'>0$ and $ C'(\varepsilon , \alpha )>0 $ we have 
         \begin{eqnarray} 
           \nonumber 
           \bigg|\frac{\partial^2 f}{\partial x^2}
           \left( t,x,u_\varepsilon(t,x),\frac{\partial u_\varepsilon}{\partial x}(t,x)\right)
           \bigg|
           & \leq & 
           C + \frac{C}{(1+x^2)^{\alpha+1} }
           \Big(\frac{\partial u_\varepsilon}{\partial x}(t,x)\Big)^2
           \\
           \label{df.xx}
            & \leq & 
           C' + C'(\varepsilon , \alpha ) e^{-2Kt} {(1+x^2)^\alpha}
         \end{eqnarray} 
         and
         \begin{eqnarray} 
\nonumber 
           \bigg|\frac{\partial^2 f}{\partial x \partial z}
           \left(t,x,u_\varepsilon(t,x),\frac{\partial u_\varepsilon}{\partial x}(t,x)\right)
           \bigg|
           & \leq &
           C +\frac{C}{(1+x^2)^{\alpha+1} }
           \Big|\frac{\partial u_\varepsilon}{\partial x}(t,x)\Big|
\\
           \label{df.xz}
& \leq &
C' + C'(\varepsilon , \alpha ) e^{-Kt}, 
        \end{eqnarray} 
$(t,x) \in [0,T]\times \real$. 
         By the relation $ v_K(t,x):=e^{-Kt}(1+x^2)^{\alpha + 1}$ and
         Condition~(\hyperlink{H5}{$H_5$}),
	we find 
	\begin{align}
	\label{a}  \bigg| \left(\frac{\partial^2 v_K}{\partial x^2}(t,x) \right)^2\frac{\partial^2 f}{\partial z^2}\Big(t,x,u_\varepsilon(t,x),\frac{\partial u_\varepsilon}{\partial x}(t,x)\Big)\bigg| 
	  & \leq C(\alpha)e^{-2Kt}(1+x^2)^{\alpha},
          	\end{align}
	for some constant $ C(\alpha) >0 $. 
         Next, we note that for some $C'(\alpha )>0$ we have 
$$
           \left| \frac{\partial^2 v_K}{\partial x^2}(t,x) \right| \leq C'(\alpha )e^{-Kt}(1+x^2)^{\alpha},
\qquad 
           \left| \frac{\partial^3 v_K}{\partial x^3}(t,x) \right|\leq C'(\alpha)e^{-Kt}(1+|x|)(1+x^2)^{\alpha -1},
           $$
           and
           $$
           \left|\frac{\partial^4 v_K}{\partial x^4}(t,x)\right|\leq C'(\alpha)e^{-Kt}(1+x^2)^{\alpha-1},
         $$
        hence from (\hyperlink{H3}{$H_3$}), \eqref{cc2} and \eqref{df.xz}
        we check that
\begin{equation}
  \nonumber 
  \left|\frac{1}{2}\frac{\partial^2 \sigma^2}{\partial x^2} + \frac{\partial f}{\partial y}
+2\frac{\partial^2 f}{\partial x \partial z}
\right|\left|\frac{\partial^2 v_K}{\partial x^2}(t,x)\right| \leq C'( \varepsilon , \alpha )e^{-Kt}(1+x^2)^{\alpha}.
\end{equation}
Similarly, thanks to the conditions~(\hyperlink{H3}{$H_3$}) and \eqref{cc2} 
\begin{equation}
  \nonumber 
  \left|\frac{\partial f}{\partial z}
+2\sigma(t,x)\frac{\partial \sigma}{\partial x}(t,x)
\right|\left|\frac{\partial^3 v_K}{\partial x^3}(t,x)\right| \leq C''(\alpha )e^{-Kt}(1+x^2)^{\alpha}, 
\end{equation}
 and 
\begin{equation}
  \label{f} \frac{1 }{2 }\sigma^2(t,x)
  \left|\frac{\partial^4 v_K}{\partial x^4}(t,x)\right|\leq C''(\alpha )e^{-Kt}(1+x^2)^{\alpha}.
\end{equation}
Combining \eqref{a}-\eqref{f}, the
inequality $\displaystyle \frac{\partial^3 v_K}{\partial x^2 \partial t}(t,x)\leq -KC'''(\alpha)e^{-Kt}(1+x^2)^{\alpha}$
for some constant $C'''(\alpha )>0$, and \eqref{df.xx}, 
 and letting $K:=1/T$ we obtain 
\begin{align*}
  S'&\ge ( e^{-KT} KC'''(\alpha ) - 2C''(\alpha ) - 2 C'( \varepsilon , \alpha ) -C(\alpha )
  ) (1+x^2)^{\alpha} - C'
\\
& \ge \left( \frac{C'''(\alpha )}{e T} -2C''( \alpha ) - 2C'( \varepsilon , \alpha )
-C(\alpha ) \right) (1+x^2)^{\alpha} - C'
\\
&\ge 0,\qquad (t,x)\in [0,T]\times \real, 
\end{align*}
 for small enough $T=T(\alpha,\varepsilon)$. 
\end{Proof}
\section{Monotonicity and continuous dependence results} 
\label{s6}
\subsubsection*{Monotonicity of FBSDEs}
In Proposition~\ref{mono1} and Proposition~\ref{nondecreasing} 
 we prove the monotonicity results needed
in the proofs of Theorem~\ref{pp4},
Corollaries~\ref{pp4.1}-\ref{pp5}, \ref{pp7}
and Theorem~\ref{pp2}.
We apply
Lemma~\ref{linear} below to derive monotonicity results for FBSDE flows of the form  
\begin{subequations}
\begin{empheq}[left=\empheqlbrace]{align}
\label{Forward2}
&  dX_s^{t,x} =\mu\big(t,X_s^{t,x}\big)ds + \sigma\big(s,X_t^{t,x}\big)dB_s, \quad X_t^{t,x}=x,
\\
\nonumber
\\
\label{Forward2.0}
 & dY_s^{t,x} =-g\big(s,X_s^{s,x},Y_s^{t,x},Z_s^{t,x}\big)ds + Z_s^{t,x}dB_s, \quad Y_T^{t,x} = \phi(X_T^{t,x}),
\end{empheq}
\end{subequations}
 $0\leq t \leq s \leq T$. 
We first prove a 
monotonicity result for the solution $(X_s^{t,x})_{s\in [t,T]}$
of the SDE~\eqref{Forward2},
which will be used to
prove non-decreasing property of $ Y_s^{t,x} $ and $  u(t,x) $
in Proposition~\ref{nondecreasing}. 
  \begin{prop}
    \label{mono1}
    Under the assumption
        (\hyperlink{A1}{$A_1$})
   the solution $\big(X_s^{t,x}\big)_{s\in [t,T]}$
  of \eqref{Forward2} 
  is a.s. non-decreasing in $x$ for all $t\in [0,T]$ and $s\in [t,T]$. 
\end{prop}
\begin{Proof}
  Let
  $ \widehat{X}^{t,x,y}_s:=   X_s^{t,y} -X_s^{t,x}  $
  for $x \leq y$, $s\in [t,T]$,
  and consider
the processes
$$
\widehat{\mu}_u:=\frac{\mu(u,X_u^{t,y} ) -\mu(u,X_u^{t,x} )}{X_u^{t,y} -X_u^{t,x}}\mathbf{1}_{\{X_u^{t,y} \ne X_u^{t,x}\}}
\ \ 
\mbox{and}
\ \ 
\widehat{\sigma}_u:=\frac{\sigma(u,X_u^{t,y} ) -\sigma(u,X_u^{t,x} )}{X_u^{t,y}-X_u^{t,x}}\mathbf{1}_{\{X_u^{t,y} \ne X_u^{t,x}\}}, 
$$
 $u\in [t,T]$.
We note that the processes
$(\widehat{\mu}_u)_{u\in [t,T]}$ and
$(\widehat{\sigma}_u)_{u\in [t,T]}$ are
bounded since $ \mu(t,x) $ and $ \sigma(t,x) $ are Lipschitz in $x$,
and that 
$(\widehat{X}^{t,x,y}_s)_{s\in [t,T]}$ satisfies the equation 
$$ \widehat{X}^{t,x,y}_s= y-x + \int_t^s\widehat{\mu}_u\widehat{X}^{t,x,y}_udu +\int_t^s\widehat{\sigma}_u\widehat{X}^{t,x,y}_udB_u,
\qquad
s\in [t,T],  $$
 which yields 
 $$ \widehat{X}^{t,x,y}_s=(y-x)\exp\Big(\int_t^s\widehat{\mu}_udu  +\int_t^s\widehat{\sigma}_udB_u -\frac{1}{2}\int_t^s\widehat{\sigma}_u^2du \Big) \ge 0, \quad 0\leq t\leq s \leq T.
 $$
\end{Proof}
\subsubsection*{Monotonicity of nonlinear PDE solutions}
The next monotonicity result is used for the proofs of Theorem~\ref{pp4},
Corollaries~\ref{pp4.1}-\ref{pp5}, \ref{pp7} and Theorem~\ref{pp2}.
\begin{prop} \label{nondecreasing}   	
  Assume that the coefficients $ \mu $, $ \sigma$, 
  $ g$ and $\phi$ satisfy (\hyperlink{A1}{$A_1$})-(\hyperlink{A4}{$A_4$}). 
  If $ \phi(x) $ and $ g(t,x,y,z) $ are non-decreasing in $ x\in\real $ for all
  $t\in[0,T]$ and $y,z \in \real $, then the solution $\big(Y_s^{t,x}\big)_{s\in [t,T]}$
	of \eqref{Forward2.0} 
	is a.s. non-decreasing in $x$ for all $s\in [t,T]$. As a consequence,
        if $u(t,x) $ is solution of the backward PDEs \eqref{du12}, 
        then $u(t,x) $ is also a non-decreasing function of $x\in \real$
        for all $t\in [0,T]$.
\end{prop}
\begin{Proof}
  	Letting 
                $\widehat{X}_s = X_s^{t,y}-X_s^{t,x}$,
                $\widehat{Y}_s = Y_s^{t,y}-Y_s^{t,x}$,
                $\widehat{Z}_s = Z_s^{t,y}-Z_s^{t,x}$ and
                $\widehat{Y}_T = \phi_2\big(X_T^{t,y}\big)-\phi_1\big(X_T^{t,x}\big)$,
                we have 
 		\begin{equation}
 		  \widehat{Y}_s = \widehat{Y}_T + \int_s^T \left(
	          g\big(u,X_u^{t,y},Y_u^{t,y},Z_u^{t,y}\right)-g\big(u,X_u^{t,x},Y_u^{t,x},Z_u^{t,x}\big)\big)
                  du-\int_s^T \widehat{Z}_udB_u. \nonumber
 		\end{equation}
 	Defining the processes $a_u, b_u$, and $c_u $ as 
 		\begin{align*} 
                  & 		a_u := \frac{g\big(u,X_u^{t,y},Y_u^{t,y},Z_u^{t,y}\big)-g\big(u,X_u^{t,x},Y_u^{t,y},Z_u^{t,y}\big)}{X_u^{t,y}-X_u^{t,x}}\mathbf{1}_{\{X_u^{t,y}\ne X_u^{t,x}\}},
                  \\
                   & 
 		  b_u := \frac{g\big(u,X_u^{t,x},Y_u^{t,y},Z_u^{t,y}\big)-g\big(u,X_u^{t,x},Y_u^{t,x},Z_u^{t,y}\big)}{Y_u^{t,y}-Y_u^{t,x}}\mathbf{1}_{\{Y_u^{t,y}\ne Y_u^{t,x}\}},
                  \\
                   & 
 		  c_u := \frac{g\big(u,X_u^{t,x},Y_u^{t,x},Z_u^{t,y}\big)-g\big(u,X_u^{t,x},Y_u^{t,x},Z_u^{t,x}\big)}{Z_u^{t,y}-Z_u^{t,x}}\mathbf{1}_{\{Z_u^{t,x}\ne Z_u^{t,x}\}},
                 \qquad
                  u\in [0,T], 
 		\end{align*} 
                which are $({\cal F}_t)_{t\in [0,T]}$-adapted and bounded
                since $g(u,x,y,z) $ is Lipschitz,
                and using the decomposition 
 		\begin{align*} 
& g\left(u,X_u^{t,y},Y_u^{t,y},Z_u^{t,y}\right)-g\big(u,X_u^{t,x},Y_u^{t,x},Z_u^{t,x}\big)
                                    \\
                   & =  g\left(u,X_u^{t,y},Y_u^{t,y},Z_u^{t,y}\right)-g\left(u,X_u^{t,x},Y_u^{t,y},Z_u^{t,y}\right) 
                  + g\left(u,X_u^{t,x},Y_u^{t,y},Z_u^{t,y}\right)-g\left(u,X_u^{t,x},Y_u^{t,x},Z_u^{t,y}\right) 
 		  \\
                   &  
 	\quad 	 + g\big(u,X_u^{t,x},Y_u^{t,x},Z_u^{t,y}\big)-g\big(u,X_u^{t,x},Y_u^{t,x},Z_u^{t,x}\big), 
 		\end{align*}
 we have 
 		\begin{equation} 
 		  \widehat{Y}_s = \widehat{Y}_T + \int_s^T \big( a_u\widehat{X}_u + b_u\widehat{Y}_u + c_u\widehat{Z}_u\big) du-\int_s^T \widehat{Z}_udB_u,
                  \qquad s\in [0,T]. \nonumber
 		\end{equation}
 		Hence, by Lemma~\ref{linear} below we get
 		\begin{equation}
	          \label{jkdlsf}
                   \widehat{Y}_s =\frac{1}{\Gamma_s} \E\left[
                  \Gamma_T\widehat{Y}_T+ \int_s^T 
                  a_u\widehat{X}_u \Gamma_u du \ \Big| \ \cF_s \right],
                \qquad
                s\in [t,T],
                \end{equation} 
 		where
 		$$\Gamma_s := \exp\left(\int_t^sc_udB_u-\frac{1}{2}\int_t^sc_u^2du+\int_0^sb_udu\right),
                \qquad
                s\in [t,T]. 
                $$
                By Proposition~\ref{mono1}
  the solution $(X_s^{t,x})_{s\in [t,T]}$
  of the forward SDE~\eqref{Forward2} satisfies
 $ \widehat{X}_s = X_s^{t,y} - X_s^{t,x} \geq 0$ for all $s \in [t,T]$ 
  if $x \leq y $, and
since $g(s,x,y,z) $ is non-decreasing in $x $
we have $ a_s \ge 0$ a.s.,
  $s\in [t,T]$. 
                Since $\phi(x) $ is non-decreasing we have
  $\widehat{Y}_T = 
                \phi\left(X_T^{t,x}\right) - \phi\left(X_T^{t,y}\right) \ge 0$ a.s.,
                hence by \eqref{jkdlsf}
                we have 
 $\widehat{Y}_s = Y_s^{t,x} - Y_s^{t,y} \geq 0$,
 $s\in [t,T]$,
                if $ x \leq y $, which implies the
                monotonicity of $ \big(Y_s^{t,x}\big)_{s\in [t,T]}$,
                therefore we also get $ u(t,x) \leq u(t,y)$, $x \leq y$,
                $t \in [0,T]$, since $ u(t,x)=Y_t^{t,x} $.  
\end{Proof}
\subsubsection*{Linear FBSDEs}
The following Lemma~\ref{linear}, which has been used in the
proof of Proposition~\ref{nondecreasing}, 
extends a classical 
 result from linear BSDEs to linear FBSDEs. 
 Let $(X_t)_{t\in [0,T]}$ satisfy the forward diffusion equation 
 \begin{equation}
   \label{Forward}
   dX_t =\mu\left(t,X_t\right)dt + \sigma\left(t,X_t\right)dB_t,
    \end{equation}
 where $\mu, \sigma$ satisfy (\hyperlink{A1}{$A_1$}),
 with associated linear backward SDE 
 \begin{equation}
   \label{linearFB}
 dY_t = -\left(a_tX_t+ b_tY_t + c_t Z_t + k_t \right)dt + Z_tdB_t,  
 \end{equation}
 with terminal condition $Y_T = \phi(X_T)$,
 where $(a_t)_{t\in [0,T]}$, $(b_t)_{t\in [0,T]}$ and $(c_t)_{t\in [0,T]}$ are
 real-valued, $({\cal F}_t)_{t\in [0,T]}$-adapted bounded processes,
 and $(k_t)_{t\in [0,T]}$ is a real-valued
 $({\cal F}_t)_{t\in [0,T]}$-adapted process such that
 $$
 \E \left[ \int_0^T k_t^2dt\right] <\infty.
 $$
 
 \begin{lemma}
   \label{linear}
   Let $(X_t)_{t\in [0,T]}$ be the solution of \eqref{Forward},
   and let $(Y_t,Z_t)_{t\in [0,T]}$ be the solution of \eqref{linearFB}. 
   Then
   the process $(Y_t)_{t\in [0,T]}$ is given in explicit form
    as 
 		\begin{equation}
               \nonumber 
 		  Y_t = \frac{1}{\Gamma_t}\E\left[ \Gamma_T\phi\left(X_T\right)+ \int_t^T \left(a_sX_s + k_s\right)\Gamma_sds \ \! \Big| \! \ {\cal F}_t\right], 
 		\end{equation}
 where $(\Gamma_t)_{t\in [0,T]}$ is the geometric Brownian motion 
 \begin{equation}
   \nonumber
   \Gamma_t := \exp\left(\int_0^tb_sds
   + \int_0^tc_sdB_s-\frac{1}{2}\int_0^tc_s^2ds
   \right), \qquad t\in [0,T]. 
 \end{equation}
 \end{lemma}
  \begin{Proof}
We have 
 	\begin{eqnarray}
 	d (\Gamma_sY_s) &=& \Gamma_sdY_s+ Y_sd\Gamma_s +  d\langle\Gamma_s,Y_s\rangle \nonumber \\
 	&=& \Gamma_s\left(-\left(a_sX_s+ b_sY_s + c_s Z_s + k_s \right)ds + Z_sdB_s\right)
        + Y_s\Gamma_s\left(b_s ds + c_sdB_s\right) + c_sZ_s\Gamma_s ds \nonumber \\
 	&=& -\left(a_sX_s\Gamma_s + k_s\Gamma_s\right)ds + \left(c_sY_s\Gamma_s+Z_s\Gamma_s\right) dB_s,
        \nonumber  
 	\end{eqnarray}
 	hence 
 	\begin{equation}
          \label{dkjld} 
 	  \Gamma_TY_T-\Gamma_tY_t = -\int_t^T \left(a_sX_s\Gamma_s + k_s\Gamma_s\right)ds + \int_t^T \left(c_sY_s\Gamma_s+Z_s\Gamma_s\right) dB_s,
 	\end{equation} 
 	and by taking conditional expectation on both sides of \eqref{dkjld}
        we find 
 	\begin{eqnarray*} 
	Y_t & = & 
        \E \left[ Y_t \mid \cF_t\right]
        \\
         & = & 
        \frac{1}{\Gamma_t}
        \E \left[ \Gamma_TY_T \mid \cF_t\right]
        + 
        \frac{1}{\Gamma_t}
        \E \left[ \int_t^T \left(a_sX_s\Gamma_s + k_s\Gamma_s\right)ds
          \ \! \Big| \ \!
          \cF_t\right]
\\
        & = &
        \frac{1}{\Gamma_t}\E \left[
          \Gamma_T\phi\left(X_T\right)+\int_t^T \left(a_sX_s + k_s\right) \Gamma_sds
          \ \! \Big| \ \!
          \cF_t\right],
        \qquad t\in [0,T].
 	\end{eqnarray*} 
\end{Proof}
\subsubsection*{Continuous dependence of FBSDE solutions} 
The next Proposition~\ref{dependence}
result extends the argument of Theorem~9.7 in \cite{mishura}
to the setting of FBSDEs. 
Other continuous dependence results are available in the literature
such as Theorem~3.3 of \cite{jakobsen}, which however 
requires uniform estimates on coefficients. 
\begin{prop} \label{dependence}
  Consider the family of forward-backward stochastic differential equations
$$ 
  \left\{
  \begin{array}{l} 
	  \nonumber
          	  \displaystyle
          X_{n,t}= X_{n,0} + \int_0^t \mu_n (s,X_{n,s} )ds + \int_0^t \sigma_n (s,X_{n,s} )dB_s, 
	 \\
         \\
	\nonumber 
	 \displaystyle
          Y_{n,t}= \phi_n(X_{n,T})+\int_{t}^{T}g_n (s,X_{n,s}, Y_{n,s}, Z_{n,s} )ds - \int_t^T Z_{n,s} dB_s, 
  \end{array}
  \right. 
$$ 
where, for every $n\geq 1$, 
the coefficients $ \mu_n $, $ \sigma_n$, 
$ g_n$ and $\phi_n$
satisfy (\hyperlink{A1}{$A_1$})-(\hyperlink{A4}{$A_4$})
for a same $C>0$. 
Assume the pointwise convergences $X_{n,0} \to X_0$ and 
\begin{equation*}
\mu_n(t,x) \to \mu(t,x), \quad
\sigma_n(t,x) \to \sigma(t,x), \quad
g_n(t,x,y,z) \to g(t,x,y,z), 
\end{equation*}
 for all $ t\in [0,T]$, and $x,y,z \in \real $ as $n \to \infty$, 
 and the strong convergence $\phi_n(x_n) \to \phi(x)$
 whenever $x_n \to x \in \real$,
 where $ \mu $, $ \sigma $, 
  $g$ and $\phi$ satisfy  
(\hyperlink{A1}{$A_1$})-(\hyperlink{A4}{$A_4$}) for a 
same constant $C>0$. Then for all $ t \in [0,T] $ we have
\begin{equation*}
	\lim\limits_{n \to \infty}\E\big[ \big|Y_{n,t}-Y_t\big|^2 \big] = 0, 
\end{equation*}
 where $(Y_t)_{t\in \real_+}$ is solution of the FBSDE
\begin{subequations}
	\begin{empheq}[left=\empheqlbrace]{align}
	  \nonumber
	  & X_t= X_0 + \int_0^t \mu (s,X_s)ds + \int_0^t \sigma (s,X_s)dB_s, 
	\\
	\nonumber 
	& Y_t= 
	\phi(X_T)+\int_t^Tg(s,X_s, Y_s, Z_s)ds - \int_t^TZ_sdB_s. 
	\end{empheq}
\end{subequations}
\end{prop} 
\begin{Proof}
  For $ n\ge 1 $, let 
  $$
  \widehat{Y}_{n,t} := Y_{n,t}-Y_t,
  \quad
  \widehat{X}_{n,t} := X_{n,t}-X_t,
  \quad
  \widehat{Z}_{n,t} := Z_{n,t}-Z_t,
  \ \ \ \! 
  \mbox{and} \quad
  \widehat{\phi}_n(x) := \phi_n(x)-\phi(x), 
  $$ 
  $t\in [0,T]$, $x\in \real$.
  Applying the It\^{o} formula to $ |\widehat{Y}_{n,t} |^2 $ and
  taking expectation on both sides yields
\begin{eqnarray*}
  \E\big[ |\widehat{Y}_{n,t}|^2 \big] & = &
  \E \big[ \big( \phi_n(X_{n,T})-\phi(X_T)\big)^2 \big]
\\
& &
+
2 \E \left[
  \int_t^T \widehat{Y}_{n,s} (
  g_n(s,X_{n,s}, Y_{n,s}, Z_{n,s})-g(s,X_s, Y_s, Z_s)
  )
  ds
  \right]
- \E \left[
  \int_t^T|\widehat{Z}_{n,s}|^2 ds
  \right]. 
\end{eqnarray*}
By the inequality $ 2 ab \leq (6C^2) a^2 + b^2/(6C^2)$, we have 
\begin{eqnarray*}
  \lefteqn{
   2 \widehat{Y}_{n,s}
  (
  g_n(s,X_{n,s}, Y_{n,s}, Z_{n,s})-g(s,X_s, Y_s, Z_s)
  )
\leq  6C^2 |\widehat{Y}_{n,s}|^2 
        }
        \\
        & &
      \ \ \ \ \   \ \ \ \ \ \ \ \ \ \ \ \ \ \ \
        +\frac{1}{6C^2} (
  g_n(s,X_{n,s}, Y_{n,s}, Z_{n,s} )
 -g(s,X_s, Y_s, Z_s) )^2.  
\end{eqnarray*}
Next, letting 
$$
\widehat{g}_n(t,x,y,z) := g_n(t,x,y,z)-g(t,x,y,z),
\qquad
t\in [0,T],
\quad
x,y,z\in \real,
$$ 
 we have 
\begin{eqnarray*}
  \lefteqn{
      (
  g_n(s,X_{n,s},Y_{n,s},Z_{n,s})-g(s,X_s,Y_s,Z_s))^2 
}
  \\
   & \leq & 2
  (g_n(s,X_{n,s},Y_{n,s},Z_{n,s})-g_n(s,X_s,Y_s,Z_s))^2 + 2 ( \widehat{g}_n(s,X_s,Y_s,Z_s) )^2
\\
&\leq & 2 C^2(|X_{n,s}-X_s|+|Y_{n,s}-Y_s|+|Z_{n,s}-Z_s|)^2 + 2 ( \widehat{g}_n(s,X_s,Y_s,Z_s) )^2
\\
&\leq & 6C^2(
|\widehat{X}_{n,s} |^2+|\widehat{Y}_{n,s} |^2+|\widehat{Z}_{n,s} |^2)
 + 2 ( \widehat{g}_n(s,X_s,Y_s,Z_s) )^2, 
\end{eqnarray*}
by the inequality $ (a+b+c)^2 \leq 3(a^2+b^2+c^2) $.
Combining the above estimates, we find 
\begin{eqnarray*}
  \E\big[
    |\widehat{Y}_{n,t}|^2
    \big]
  &\leq & \E \big[ \big( \phi_n(X_{n,T})-\phi(X_T)\big)^2 \big]
  \\
   & & 
  + (6C^2+1) \E \left[
    \int_t^T |\widehat{Y}_{n,s}|^2 ds \right]
  + \E \left[ \int_t^T |\widehat{X}_{n,s}|^2 ds \right]
  \\
  & & + \frac{1}{3C^2} \E \left[ \int_t^T 
    \big( \widehat{g}_n(s,X_s,Y_s,Z_s) \big)^2 ds \right], 
\end{eqnarray*}
 which yields 
\begin{align*}
  \E\big[ |\widehat{Y}_{n,t}|^2 \big]
  &\leq C''\left(
  \E \big[ \big( \phi_n(X_{n,T})-\phi(X_T)\big)^2 \big]
  + \int_t^T \E \big[ |\widehat{X}_{n,s}|^2 \big] ds  
  + \E \left[ \int_t^T \big( \widehat{g}_n(s,X_s,Y_s,Z_s) \big)^2 ds \right] \right)
\end{align*}
by Gronwall's inequality, and therefore 
\begin{align*}
 \E\big[
    |\widehat{Y}_{n,t}|^2
    \big]
 &\leq C''
 \E \big[ \big( \phi_n(X_{n,T})-\phi(X_T)\big)^2 \big]
  \\
  & \quad +  C'' (T-t) \sup_{s\in [t,T]} \E \big[ |\widehat{X}_{n,s}|^2 \big]+
  C'' \E \left[ \int_t^T \big( \widehat{g}_n(s,X_s,Y_s,Z_s) \big)^2 \big] ds \right]. 
\end{align*}
We note that since $X_{n,0} \to X_0$,
$\mu_n(t,x) \to \mu $ and $ \sigma_n(t,x) \to \sigma(t,x)$ pointwise when
$ n \to \infty $, by Theorem~9.7 of \cite{mishura} 
we have $\lim \limits_{n \to \infty} \E \big[ \sup_{t\in [0,T]} |\widehat{X}_{n,t}|^2 \big] =0$. 
Hence, by the condition $ |\widehat{g}_n(t,x,y,z)|\leq 2C(|x|+|y|+|z|)$
and the pointwise limit 
$\lim \limits_{n \to \infty} \widehat{g}_n(t,x,y,z) = 0  $,
by Lebesgue dominated convergence we find 
\begin{equation*}
  \lim\limits_{n \to \infty}
  \E \left[ \int_t^T \big( \widehat{g}_n(s,X_s,Y_s,Z_s) \big)^2 \big] ds
  \right] = 0, \quad t\in[0,T].
\end{equation*}
Finally, by the strong convergence of $(\phi_n)_{n\geq 1}$ to $\phi$ 
and the uniform integrability
$$\sup_{n\geq 1} \E \bigg[ \sup_{t\in [0,T]} |X_{n,t}|^{2p} \bigg] < \infty,
\qquad
p \geq 1,
$$
see Theorem~9.2 in \cite{mishura}, we obtain
 \begin{equation*}
  \lim\limits_{n \to \infty}
 \E \big[ \big( \phi_n(X_{n,T})-\phi(X_T)\big)^2 \big] = 0, 
\end{equation*}
 and we conclude to 
\begin{equation*}
\lim\limits_{n \to \infty}\E\left[ \left|Y_{n,t}-Y_t\right|^2 \right] = 0, \quad t\in[0,T].
\end{equation*}
\end{Proof}

\footnotesize

\def\cprime{$'$} \def\polhk#1{\setbox0=\hbox{#1}{\ooalign{\hidewidth
  \lower1.5ex\hbox{`}\hidewidth\crcr\unhbox0}}}
  \def\polhk#1{\setbox0=\hbox{#1}{\ooalign{\hidewidth
  \lower1.5ex\hbox{`}\hidewidth\crcr\unhbox0}}} \def\cprime{$'$}

\end{document}